\documentclass[]{amsart}

\usepackage{amssymb}
\usepackage{amsmath}
\usepackage{amsthm}
\usepackage{amscd}

\usepackage{overpic}

\long\def\comment#1\endcomment{}
\usepackage{color}

\newcommand {\cals}{\mathcal S} 

\newcommand {\iv}{^{-1}}
\newcommand{\Int}{{\mathrm{Int}}\, }
\newcommand{\thin}{{\mathrm{Folded}}}

\newcommand{\dv}{{\mathrm{div}}}

\newcommand{\Dv}{{\mathrm{Div}}}
\newcommand{\Is}{{\mathrm{Iso}}}

\newcommand{\Rd}{{\mathrm{Rad}}}
\newcommand{\vol}{{\mathrm{Vol}}}
\newcommand{\fvol}{{\mathrm{FillVol}}}
\newcommand{\frad}{{\mathrm{FillRad}}}
\newcommand{\Rk}{{\mathrm{Rank}}}

\newcommand{\vf}{\varphi}
\newcommand {\nn}{\mathcal N} 
\newcommand{\onn}{\overline{{\mathcal N}}}
\newcommand {\calm}{\mathcal M} 

\def\<{\left\langle}
\def\>{\right\rangle}

\newcommand {\rgot}{{\mathfrak r}}
\newcommand {\sgot}{{\mathfrak s}}

\newcommand {\cgot}{{\mathfrak c}}
\newcommand {\bgot}{{\mathfrak b}}
\newcommand {\dgot}{{\mathfrak d}}

\newcommand{\me}{\medskip}

\newcommand {\q}{\mathfrak q} 

\newcommand{\fh}{\mathfrak{h}}

\newcommand{\fk}{\mathfrak{k}}

\newtheorem{thm}{Theorem}[section]
\newtheorem{prop}[thm]{Proposition}
\newtheorem{cor}[thm]{Corollary}

\newtheorem{lem}[thm]{Lemma}
\newtheorem{cvn}[thm]{Convention}

\newtheorem{notation}[thm]{Notation}

\theoremstyle{definition} \newtheorem{defn}[thm]{Definition}

\theoremstyle{remark}
\newtheorem{rmk}[thm]{Remark}

\def\square{\hfill${\vcenter{\vbox{\hrule height.4pt \hbox{\vrule
width.4pt height7pt \kern7pt \vrule width.4pt} \hrule
height.4pt}}}$}

\newcommand{\tsh}[1]{\left\{\kern-.9ex\left\{#1\right\}\kern-.9ex\right\}}

\def\R{{\mathbb R}}

\def\N{{\mathbb N}}

\def\dd{{\mathcal D}}

\def\calc{{\mathcal C}}

\def\V{{\mathcal V}}

\def\g{{\mathcal G}}

\def\calr{{\mathcal R}}

\def\nn{{\mathcal N}}

\def\fh{{\mathfrak h}}

\def\diam{{\rm{diam}}}

\def\dist{{\rm{dist}}}

\def\<{\langle}
\def\>{\rangle}

\long\def\Restate#1#2#3#4{
\medskip\par\noindent
{\bf #1 \ref{#2} #3} {\it #4}\par\medskip }

\definecolor{darkgreen}{cmyk}{1,0,1,.2}

\def\lam{\lambda}

\newcommand\intersect\cap
\newcommand\infinity\infty
\newcommand\wt\widetilde

\newcommand\inject\hookrightarrow
\newcommand\union\cup

\newcommand{\co}{\colon\thinspace}

\newcommand\join\Lambda
\newcommand\cross\times

\newcommand\lub\vee
\newcommand\glb\wedge

\renewcommand\paragraph[1]{\medskip\textbf{#1} }

\begin{document}
\title[Combinatorial isoperimetry and divergence] 
{Combinatorial higher dimensional 
isoperimetry and divergence}
\author{Jason Behrstock}
\thanks{The research of the first author was supported
as an Alfred P. Sloan Fellow, a Simons Fellow, and by
NSF grant \#DMS-1006219}
\address{Lehman College and the Graduate Center,
City University of New York,
U.S.A.}
\curraddr{Barnard College, Columbia University, New York, New York,
USA} 
\author{Cornelia Dru\c{t}u}\thanks{The research of the second author was supported in part by
the EPSRC grant ``Geometric and analytic aspects of infinite groups", by the project ANR Blanc ANR-10-BLAN 0116, acronym GGAA, and by the LABEX CEMPI}

\address{Mathematical Institute,
Radcliffe Observatory Quarter,
Oxford OX2 6GG, U. K.}


\date{\today}

\begin{abstract} 
In this paper we provide a framework for the study of 
isoperimetric problems in finitely generated group, through a 
combinatorial study of universal covers of compact
simplicial complexes.  We show that, when estimating filling
functions, one can restrict to simplicial spheres of particular
shapes, called ``round'' and ``unfolded'', provided that a bounded
quasi-geodesic combing exists.  We prove that the problem of
estimating higher dimensional divergence as well can be restricted to
round spheres. Applications of these results include a combinatorial 
analogy of the Federer--Fleming 
inequality for finitely 
generated groups, the construction of examples of $CAT(0)$--groups 
with higher dimensional divergence equivalent to $x^d$ for every degree $d$ \cite{BehrstockDrutu:HigherFilling2}, and a proof of the fact that for bi-combable groups the
filling function above the quasi-flat rank is asymptotically linear \cite{BehrstockDrutu:Linear}.
\end{abstract}

\maketitle

\section{Introduction}

The \emph{$k$--dimensional isoperimetric (or filling) function} of a
space $X$, denoted in this paper $\Is_k$, measures the smallest volume
of a $(k+1)$--dimensional ball needed to fill a $k$--dimensional
sphere of a given area.  There is a whole range of \emph{filling
functions}, from the ones mentioned above, to filling functions of the
form $\Is_V$, measuring how copies of a given manifold $\partial V$
are filled by copies of a manifold $V$, where $V$ is an arbitrary
$(k+1)$--dimensional connected compact sub-manifold with boundary of
$\R^{k+1}\, $.  The notion of ``volume'' that is being used also
varies.

A particularly significant type of filling function, especially in the
presence of non-positive curvature, are the \emph{divergence
functions}. These functions measure the volume of a filling of a 
sphere, where the filling is required to avoid a ball of large size. 
These functions, in some sense, describe the
spread of geodesics, and the filling near the boundary at infinity.

Traditionally, the topic of filling (of spheres, hypersurfaces, cycles
etc) belongs to Riemannian geometry and Geometric measure theory.
More recently, it has made its way in the study of infinite groups, in
which the most appropriate framework is that of simplicial complexes
and simplicial maps; this setting arises naturally in the context of groups whose
Eilenberg-Mac Lane space has a finite $(k+1)$--skeleton.  In a 
simplicial setting, analytic arguments and tools are no longer
available and analogous tools must be established anew.

A main difficulty in obtaining estimates for higher dimensional
filling functions, whether in the context of Riemannian geometry, 
Geometric measure theory, or elsewhere,
comes from the fact that, unlike in the 1-dimensional case, the
knowledge of the volume of a hypersurface to fill does not bring with
it any knowledge of its shape, or even, in particular, its diameter.  This
difficulty does not occur in the case of one-dimensional filling,
hence in group theory significant results have already been found in
that case.  Isoperimetric functions have been used to characterize a
number of interesting classes of finitely presented groups.  For
instance, a group is Gromov hyperbolic if and only if its
$1$--dimensional isoperimetric function is subquadratic, and this
occurs if and only if this function is linear
\cite{Gromov:hyperbolic,Olshanskii:hip,Papasoglu:hip}.  Also, it has
been shown that automatic groups have $k$--dimensional isoperimetric
functions that are at most polynomial for every $k$, and that in the
particular case when $k=1$
their isoperimetric functions are always
at most quadratic \cite{ECHLPT}. Recently there has been important progress on
isoperimetric functions for lattices in higher rank semisimple groups, see
\cite{BestvinaEskinWortman, Leuzinger:optimalDehn, Young:slnz}.

\medskip

In this paper we prove that, under certain conditions, the study of
the minimal filling of simplicial spheres by simplicial balls can be
restricted to classes of spheres with particular shapes, which we call \emph{round} and \emph{unfolded}.  A $k$--dimensional
hypersurface $\fh$ is called $\eta$--round, for some $\eta >0$, if its
diameter is at most $ \eta \vol (\fh )^{\frac{1}{k}}\, $.  We defer 
the precise definition of an unfolded hypersurface to Definition
\ref{def:folded}, but, roughly speaking, a $k$--dimensional hypersurface
is unfolded at scale $\rho$ if within distance at most $\rho$
of every point the hypersurface looks like a $k$--dimensional disk
(as opposed to, say, looking like a $k$--dimensional cylinder that is
long and thin).

\emph{Throughout this paper we work in the context of an 
$n+1$--dimensional simplicial
complex $X$ which is $n$--connected, and which is the
universal cover of a compact simplicial complex.}

We prove that an arbitrary sphere of dimension $k\leq n$ in $X$ has a
partition into spheres with particular shapes, such that the sum of
the volumes of the spheres in the partition is bounded by a multiple
of the volume of the initial sphere.  Here a \emph{partition} 
consists of the following: finitely many spheres (more generally,
hypersurfaces) $\fh_1,\dots, \fh_n$ compose a partition of a sphere
$\fh$ if by filling all of them one obtains a ball filling $\fh$ (see
Definition \ref{def:partition}).  The hypersurfaces $\fh_1,\dots,
\fh_n$ are called \emph{contours of the partition}.

We prove the following:

\Restate{Theorem}{thm:rounddec0}{}
{Consider an integer $2\leq k\leq n$. Assume that $X$ satisfies an 
isoperimetric inequality at most Euclidean for $k-1$,
if $k\neq 3 , $ or an inequality of the form $\Is_V (x)\leq B
x^{3}$ for every compact closed surface bounding a handlebody 
$V$, if $k=2$ (where $B$ is independent of $V$).

Then for every $\varepsilon >0$ there exists a constant $\eta >0$ such that every $k$--dimensional sphere $\fh$ has a
partition with contours $\fh_1,...\fh_n$ that are $\eta$--round hypersurfaces, and contours $\rgot_1,\dots ,\rgot_m$ that are hypersurfaces of volume and filling volume zero such that
\begin{enumerate}

\item $\sum_{i=1}^n \vol (\fh_i ) \leq 2\cdot 6^{k+1} \vol (\fh )\, .$

\medskip

\item $\fh_1,...\fh_n$ and $\rgot_1,\dots ,\rgot_m$ are  contained in the tubular neighborhood $\nn_R (\fh )$, where $R= \varepsilon \vol (\fh )^{1/k}$.
\end{enumerate}

For $k=2$ all the hypersurfaces $\fh_1,...\fh_n$ and $\rgot_1,\dots ,\rgot_m$ are spheres.  

If moreover $X$ has a bounded quasi-geodesic combing then for every $k\geq 3$ as well, $\fh_1,...\fh_n$ are $\kappa \eta$--round spheres, where $\kappa$ depends only on the constants of the quasi-geodesic combing. 
}

An immediate consequence of the above is that a simplicial complex
with a bounded quasi-geodesic combing always satisfies an
isoperimetric inequality which is at most Euclidean. Further, an estimate for
the filling of round spheres yields one for arbitrary spheres.

The existence of bounded quasi-geodesic combing implies cone-type
inequalities in the sense of \eqref{eq:cone} in a straightforward way, see
e.g., \cite[Section~10]{ECHLPT} and Lemma~\ref{lem:combing}; in
dimension 1 these immediately imply a Euclidean filling function, but
applying these inequalities in higher dimensions to obtain Euclidean
filling functions requires significantly more elaborate arguments,
which we carry out here.  Along
the same lines, Gromov previously studied homological filling
functions by Lipschitz chains. In this context, Gromov proved that
such homological fillings are at least Euclidean for Hadamard
manifolds, and,
more generally, for complete Riemannian manifolds satisfying a cone-type
inequality (see for instance \eqref{eq:cone}) and for Banach spaces
\cite{Gromov(1983)}.  Gromov's results were extended to Hadamard spaces in
\cite{Wenger:isoperim}, and then to complete metric spaces satisfying
cone-type inequalities in \cite{Wenger:isoperim} via another
homological version of filling functions using integral currents.

From our result on the partition into round hypersurfaces previously described, we can deduce the following.

\Restate{Theorem}{thm:FedererFlemingGroups}{}
{
Assume that $X$ has a bounded quasi-geodesic combing.  
\begin{enumerate}
\item{\rm(Federer-Fleming inequality for
groups).} For every $k\geq 1$, the isoperimetric inequality  of $X$ is 
at most Euclidean.  For $k=2$, moreover, the supremum of filling functions 
$\Is_V (x)$ modelled over all handlebodies $V$ is at most $B
x^{3}$, where $B$ is independent of $V$. 

\medskip

\item Assume that for some $k\geq 2$ it is known that every round $k$--sphere of volume at most $Ax^{k}$ has filling
volume at most $B x^{\alpha}$ with $\alpha \in [k, k+1)$. Then $\Is_k (x)
\leq \kappa B x^{\alpha},$ where $\kappa $ is a universal constant.
\end{enumerate}
}

We were informed by S. Wenger that the Federer-Fleming inequality for combable groups can also be deduced from \cite{Wenger:isoperim} and \cite[Theorem 1, p.  435]{White}. 
%
%

The above result generalizes to show that also in the study of higher 
dimensional divergence one may
restrict the study to round spheres.  For a definition of higher
dimensional divergence we refer to Definition \ref{def:kdiv}.  We
prove the following.

\Restate{Theorem}{thm:DivRound}{}
{
Assume that $X$ is a simplicial complex of dimension $n$ endowed with a bounded quasi-geodesic combing.
For every $\varepsilon >0$ there exists $\eta >0$ such that the following holds. Consider the restricted
divergence function $\Dv^r_k(x, \delta )$, obtained by taking the supremum only over $k$--dimensional spheres that are
$\eta$--round, of volume at most $2A x^{k}$ and situated outside balls of radius $x$. 

Assume that $\Dv^r_k(x, \delta ) \leq Br^\beta$ for some $\beta \geq k+1$ and $B>0$ universal
constant.  Then the general divergence function $\Dv_k (x, \delta (1- \varepsilon ) )$ is at most $B'r^\beta$ for some $B'>0$ depending on $B, \varepsilon , \eta $ and $X$.
}

In fact we prove a more general version of the theorem above (see Section \ref{section:div}).

Theorem \ref{thm:DivRound} is a powerful tool which we expect will be 
widely used. Indeed, in \cite{BehrstockDrutu:HigherFilling2} the 
authors use this to construct CAT(0) groups with higher dimensional divergence exactly a polynomial function of our choice. 

When a bounded quasi-geodesic combing exists, another useful 
restriction to impose is to restrict one's study to fillings of spheres that
are both round and unfolded.  Note that the condition of roundedness
only forbids that a sphere stretches too much towards infinity, but it
does not guarantee that the sphere does not contain many long and thin
``fingers'' (which may eventually be spiralling, so that their diameter
satisfies the condition imposed by roundedness).  It is the condition
of ``unfoldedness'' that requires from a $k$-dimensional sphere to be
shaped like an Euclidean $k$-dimensional sphere.
In such a setting we prove the following result.

\Restate{Theorem}{thm:roundThick}{}
{
Let $X$ be a simplicial complex with a bounded quasi-geodesic combing.

\begin{enumerate}
\item Let $k\geq 2$ be an integer. If every $k$--dimensional sphere of volume at most $Ax^{k}$ that is
$\eta$--round and $\varepsilon$---unfolded at scale $\delta x$, in the sense of Definition \ref{def:folded}, has filling volume at most $B x^{\alpha}$ with $\alpha \geq k$, then $\Is_k (x) \leq A x^{\alpha},$ where $A =A (\eta , \varepsilon , \delta )$.

\item If every (closed) surface of volume at most $Ax^{2}$ that is
$\eta$--round and $\varepsilon$---unfolded at scale $\delta x$ has filling volume at most $B x^{\alpha}$ with $\alpha \geq 2$ and $B$ independent of the genus, then $\Is_V (x) \leq A x^{\alpha},$ for every handlebody $V$, where $A =A (\eta , \varepsilon , \delta )$.

\end{enumerate}
}

Unlike the case of reduction to round spheres, Theorem
\ref{thm:roundThick} is not based on a partition as in Theorem
\ref{thm:FedererFlemingGroups}. Instead, for this result, for an arbitrary sphere, 
we produce a partition with a uniformly bounded number of contours that
are spheres round and unfolded and where the remaining contours are spheres
whose respective area is at most $\varepsilon$ times the area of the
initial sphere, with $\varepsilon >0$ small.

Theorem \ref{thm:roundThick} is a major tool in our argument in
\cite{BehrstockDrutu:Linear} proving that for a bi-combable group the
filling function above the quasi-flat rank is asymptotically linear. 
There, we make use of the fact that filling spheres 
that are both 
round and unfolded descend to spheres which retain a nice geometric 
structure when passing to asymptotic cones and thereby allowing for 
clean relationships between these fillings and the rank.

Partial results in the direction of a decomposition of spheres into
special types of spheres that yield good limits in the asymptotic
cones have been obtained previously by P. Papasoglu in
\cite{Papasoglu:quadratic, Papasoglu:iso}, T. Riley in
\cite{TRiley:higherIso} and S. Wenger in \cite{Wenger:asrk}.

The plan of the paper is as follows.  In Section \ref{sect:prelim0} we
recall some basic notions and establish notation which we will use in
the paper.  In Section \ref{sect:higherIso} we recall a few facts
about filling functions and we prove an estimate of the filling
radius in terms of the filling function.  Section
\ref{sec:roundunfold} is devoted to the proof that in a simplicial
complex with bounded combing the study of the filling function can be
reduced to round and unfolded spheres.  In Section \ref{section:div},
we prove that in the presence of a combing, the study of divergence
can likewise be restricted to round spheres.

\subsection*{Acknowledgements}
The authors thank Bruce Kleiner, Enrico Leuzinger and Stefan Wenger
for useful comments on an earlier version of this paper. 
J.B.\ thanks the Mathematics Departments of Barnard/Columbia and  
Exeter College/Oxford University for their 
hospitality during the writing of this paper.

\section{Preliminaries}\label{sect:prelim0}

\subsection{General terminology and notation}\label{sect:prelim}

We begin with standard notions and notation used in the study of
quasi-isometry invariants. Consider a constant $C\geq 1$ and an integer $k\geq 1$. Given two functions $f,g$ which both map $\R_+$ to itself, we write $f\preceq_{C,k} g$ if
$$f(x) \le Cg(Cx+C)+Cx^{k}+C\; \; \mbox{ for all }x\in \R_+.
$$
 
We write $f\asymp_{C,k} g$ if and only if
$f\preceq_{C,k} g$ and $g\preceq_{C,k} f$. 
Two functions $\R_+\to \R_+$ are
said to be $k$--\emph{asymptotically equal} if there exists $C\geq 1$ s.t.
$f\asymp_{C,k} g$.  This is an equivalence relation.  

When at least one of the two functions $f,g$ involved in the relations above is an $n$--dimensional isoperimetric or divergence function, we automatically consider only relations where $k=n$, therefore $k$ will no longer appear in the subscript of the relation. When irrelevant, we
do not mention the constant $C$ either and likewise remove the corresponding
subscript.

Given $f$ and $g$ real-valued functions of one real variable,
we write $f= O(g)$ to mean that there exists a constant $L>0$ such
that $f(x) \leq L g(x)$ for every $x$; in particular $f= O(1)$ means
that $f$ is bounded, and $f= g +O(1)$ means that $f-g$ is bounded.
The notation $f= o(g)$ means that $\lim_{x\to \infty }
\frac{f(x)}{g(x)}=0\, $.

\me

In a metric space $(X, \dist )$, the {\em open $R$--neighborhood} of
a subset $A$, i.e. $\{x\in X:
\dist (x, A)<R\}$, is denoted by $\nn_R(A)$.  In particular, if $A=\{a\}$ then $\nn_R(A)=B(a,R)$
is the open $R$--ball centered at $a$. We use the notation $\onn_R (A)$ and $\bar{B}(a,R)$ to designate the
corresponding {\em closed neighborhood} and {\em closed ball}
defined by non-strict inequalities. We make the convention that $B(a,R)$ and $\bar{B}(a,R)$ are the empty
set for $R<0$ and any $a\in X$.

\me

Fix two constants $L\geq 1$ and $C\geq 0$. A map $\q\co Y\to X$ is
said to be
\begin{itemize}
  \item  $(L,C)$--\emph{quasi-Lipschitz} if
$$
\dist (\q (y),\q (y'))\leq L\dist (y,y')+C,
\hbox{ for all } y,y' \in Y\, ;
$$
  \item  an $(L,C)$--\emph{quasi-isometric embedding} if moreover
  $$
\dist (\q (y),\q (y')) \geq \frac{1}{L}(y,y')-C\hbox{ for all } y,y' \in Y\, ;
$$
  \item an $(L,C)$-\emph{quasi-isometry} if it is an $(L,C)$--quasi-isometric
embedding $\q\co Y\to X$ satisfying the additional assumption that $X\subset \nn_C(\q
(Y))$.
  \item an $(L,C)$--{\em quasi-geodesic} if it is an $(L,C)$--quasi-isometric
embedding defined on an interval of the real
line;

\item a \emph{bi-infinite $(L,C)$--quasi-geodesic} when defined on the entire real line.
\end{itemize}

In the last two cases the terminology is extended to the image of $\q$. When the constants $L,C$ are
irrelevant they are not mentioned.

We call $(L,0)$--quasi-isometries (quasi-geodesics)
$L$--\emph{bi-Lipschitz maps (paths)}.
If an $(L,C)$--quasi-geodesic $\q$ is $L$--Lipschitz then $\q$ is called an {\em $(L,C)$--almost geodesic}.
Every
$(L,C)$-quasi-geodesic in a geodesic metric space is at bounded (in
terms of $L, C$) distance from an $(L+C,C)$--almost geodesic with the
same end points, see e.g. \cite[Proposition 8.3.4]{Buragos:course}.
Therefore, without loss of generality, we assume in
this text that all quasi-geodesics are in fact almost geodesics, in
particular that they are continuous.

\me

Given two subsets $A,B\subset \R$,
a map $f\co A\to B$ is said to be
\emph{coarsely increasing} if there exists a constant $D$ such that
for each $a,b$ in $A$ satisfying $a+D<b$, we have that $f(a)\leq f(b)$. Similarly, we define
\emph{coarsely decreasing} and \emph{coarsely monotonic} maps. A
map between quasi-geodesics is coarsely monotonic if it defines a
coarsely monotonic map between suitable nets in their domain.

\me

A metric space is called
\begin{itemize}
  \item \emph{proper} if all its closed balls are compact;
  \item \emph{cocompact} if there exists a compact subset $K$ in $X$ such that all the translations of $K$ by isometries of $X$ cover $X$;
  \item \emph{periodic} if it is geodesic
and for fixed constants $L\geq 1$ and $C\geq 0$ the image of
some fixed ball under $(L,C)$--quasi-isometries of $X$ covers $X$;
  \item a \emph{Hadamard space} if $X$ is geodesic, complete, simply connected and satisfies the CAT(0)
condition;
  \item  a \emph{Hadamard manifold} if moreover $X$ is a
smooth Riemannian manifold.
\end{itemize}

\subsection{Combinatorial terminology}\label{subsec:combi}

The usual setting for defining an $n$--dimen\-sio\-nal filling
function is
that of an $n$--connected space $X$ of dimension $n+1$;
of particular interest is when
$X$ is the universal cover of a compact CW-complex $K$, with fundamental group
$G$.  By the Simplicial Approximation Theorem,
cf. \cite[Theorem 2.C.2]{Hatcher:book}, $K$ is homotopy equivalent to
a finite simplicial complex $K'$ of the same dimension; hence we may
assume that both $X$ and $K$ are simplicial.

In this paper we use the standard terminology related to simplicial complexes as it appears in \cite{Hatcher:book}. In the setting of isoperimetry problems, this terminology is used as such in \cite{Papasoglu:iso}, and it is used in a slightly more general but
equivalent form (i.e. the cells need not be simplices but rather polyhedra with
a uniformly bounded number of faces) in
\cite[p. 153]{Bridson-Haefliger}, \cite{Bridson:PolynDehn}, and
\cite[$\S 2.3$]{TRiley:higherIso}.  Note that when we speak of
simplicial complexes in what follows, we always mean their topological realisation. Throughout the paper, we assume that all simplicial complexes are connected.

 Given an $n$--dimensional simplicial complex $\calc$, we call
  the closed simplices of dimension $n$ the
  \emph{chambers} of $\calc$.  A
 \emph{gallery} in $\calc$ is a finite sequence of chambers such that two consecutive chambers share a face of dimension $n-1$.

Given a simplicial map $f\co X\to Y$, where $X,Y$ are
simplicial complexes, $X$ of dimension $n$, we call $f$--\emph{non-collapsed chambers in $X$} the chambers whose images by $f$ stay of dimension $n$. We denote by $X_\vol$ the set of $f$--non-collapsed chambers. We define the \emph{volume of }$f$ to be the (possibly infinite) cardinality of $X_\vol$.

\me

Recall that a group $G$ is \emph{of type} $\mathcal{F}_k$ if it admits
an Eilenberg-MacLane space $K(G,1)$ whose $k$-skeleton is finite.

\begin{prop}[\cite{AWP}, Proposition 2]\label{prop:typefn}
If a group $G$ acts cellularly on a CW-complex $X$, with finite
stabilizers of points and such that $X^{(1)}/G$ is finite then $G$ is
finitely generated and quasi-isometric to $X$.  Moreover, if $X$ is
$n$-connected and $X^{(n+1)}/G$ is finite then $G$ is of type
$\mathcal{F}_{n+1}$.
\end{prop}

Conversely, it is easily seen that for a group of type
$\mathcal{F}_{n+1}$ one can define an $(n+1)$-dimensional
$n$-connected simplicial complex $X$ on which $G$ acts
properly discontinuously by simplicial isomorphisms, with trivial
stabilizers of vertices, such that $X/G$ has finitely many cells.  Any
two such complexes $X,Y$ are quasi-isometric, and the quasi-isometry,
which can initially be seen as a bi-Lipschitz map between two subsets
of vertices, can be easily extended to a simplicial map $X\to Y\, $
\cite[Lemma 12]{AWP}.

A group is \emph{of type} $\mathcal{F}_\infty$ if and only if it is of
type $\mathcal{F}_k$ for every $k\in \N\, $.
It was proven in \cite[Theorem 10.2.6]{ECHLPT} that every combable group is
of type $\mathcal{F}_\infty$.

\section{Higher dimensional isoperimetric functions}\label{sect:higherIso}

\subsection{Definitions and properties}\label{sect:fillingFct}

There exist several versions of filling functions, measuring how
spheres can be filled with balls or, given a manifold pair $(M,
\partial M )$, how a copy of $\partial M$ can be filled
with a copy of $M$, or how a cycle can be filled with a chain.  The
meaning of `sphere', `manifold' or `cycle' also
varies, from the measure theoretical notion of integral current
\cite{AmbrosioKirchheim,Wenger:div,Wenger:isoperim} to that of
(singular) cellular map \cite{Bridson:PolynDehn,
TRiley:higherIso, BBFS} or of Lipschitz map defined on the proper
geometric object. For a comparison between the various versions of filling functions we refer to \cite{Groft:thesis, Groft1, Groft2}.

In the setting of finitely generated groups the most frequently used approach is to refer to a proper cocompact
action of the group on a CW-complex.  More precisely, the $n$--th
dimensional filling function is defined for groups that are \emph{of
type} $F_{n+1}$, that is groups having a classifying space with
finite $(n+1)$-skeleton.  One can define the $n$--th dimensional
filling function using the $(n+1)$-skeleton of the classifying space, or any other
$(n+1)$-dimensional complex on which the group $G$ acts properly
discontinuously cocompactly. This is due to the quasi-isometry
invariance of filling functions proved in \cite{AWP}.  Since a finite
$(n+1)$--presentation of a group composed only of simplices can always
be found, it suffices to restrict to simplicial complexes. In what follows we therefore define filling functions for simplicial complexes with a cocompact action.

A simplicial complex $X$ may be endowed with a ``large scale metric
structure'' by assuming that all edges have length one and taking the
shortest path metric on the $1$-skeleton $X^{(1)}$.  We say that a
metric space $Y$ (or, another simplicial complex~$Z$) is
$(L,C)$--\emph{quasi-isometric to} $X$ if $Y$ (respectively $Z^{(1)}$)
is $(L,C)$--quasi-isometric to~$X^{(1)}$.

\begin{cvn}
For the rest of the section, we fix a simplicial complex $X$ in which
the filling problem is to be considered.  We assume that $X$ is the
universal cover of a compact simplicial complex $K$ with
fundamental group $G$, that it has dimension $n+1$ and it is
$n$--connected. We will consider fillings in $X$ up to dimension $n$.
\end{cvn}

\begin{cvn}\label{cvn:V}
Throughout the paper, when we speak of manifolds we always mean
\emph{manifolds with a simplicial-complex structure}.

We denote by $V$ an arbitrary $m$--dimensional
connected compact sub-manifold of $\R^{m}\, $, where $m\geq 2$ is an
integer and $V$ is smooth or piecewise linear, and with boundary.  We
denote its boundary, by $\partial V$. Unless otherwise stated, the standing assumption is that $\partial V$ is connected. We denote its interior by $\Int (V)$.
\end{cvn}

Given $V$ as above, a $V$--\emph{domain in $X$} is a simplicial map $\frak d$ of $\dd$ to $X^{(m)}$, where $\dd$ is a
simplicial structure on $V$.  When the manifold $V$ is irrelevant we
simply call $\frak d$ a \emph{domain of dimension $m$} (somewhat
abusively, since it might have its entire image inside $X^{(m-1)}$);
we also abuse notation by using $\frak d$ to denote both the map and its image.

A $\partial V$--\emph{hypersurface in $X$} is a simplicial map  $\fh$ of
$\calm$ to $X^{(m-1)}$, where $\calm$ is a simplicial structure of the boundary
  $\partial V$.  Again, we abuse notation by letting $\fh$ also denote the image of the above
map, and we also call both $\fh$ and its image a \emph{hypersurface of dimension $m-1$}.

According to the terminology introduced in Section \ref{subsec:combi}, $\dd_\vol$, respectively $\calm_\vol$, is the set of $\frak d$--non-collapsed chambers (respectively $\fh$--non-collapsed chambers). The \emph{volume} of $\frak d$ (respectively $\fh$) is the cardinality of $\dd_\vol$, respectively $\calm_\vol$. Given a vertex $v$ we write $v\in \dd_\vol$ (or $v\in \calm_\vol$) to
signify that $v$ is a vertex in a non-collapsed chamber.

We sometimes say that the domain $\frak d$ is \emph{modelled on} $V$, or that it is a $V$--{\emph{domain}}, and $\fh$ is \emph{modelled on} $\partial V$, or that it is a $\partial V$--{\emph{hypersurface}}. When $V$ is a closed ball in $\R^m$, we call $\frak d$ an $m$--\emph{dimensional ball} and $\fh$ an $(m-1)$--\emph{dimensional sphere}.

We say that a domain $\frak d$ \emph{fills} a hypersurface $\fh$ if
this pair corresponds to a $(k+1)$--dimensional connected compact
smooth sub-manifold with boundary $V$ in $\R^{k+1}\, $ satisfying $\dd \cap
\partial V = \calm$ and $\frak d|_\calm = \fh$,
possibly after pre-composing $\fh$ with a simplicial equivalence of
$\calm$.

The \emph{filling volume} of the hypersurface $\fh$,
$\fvol (\fh )$, is the minimum of all the volumes of
domains filling $\fh$.  If no domain filling $\fh$ exists then we set
$\fvol (\fh ) = \infty$.

\begin{rmk}
These notions are equivalent to the ones defined in \cite{BBFS,young:nonrecursiveFV3} using admissible maps, as well as to the ones in \cite{AWP}, and those in \cite[p.153]{Bridson-Haefliger}, \cite{Bridson:PolynDehn}, \cite[$\S
2.3$]{TRiley:higherIso} that are using more polyhedra than just simplices.

Indeed, every domain and hypersurface as above is an admissible map
with the same volume.  Conversely, consider an admissible map $f\co W \to
X$ defined on an $m$-dimensional domain or boundary of a domain, i.e.
a continuous map $f\co  W \to X^{(m)}$, such that $f\iv(X^{(m)} \setminus
X^{(m-1)})$ is a disjoint union $\bigsqcup_{i\in I} B_i$ of open
$m$-dimensional balls, each mapped by $f$ homeomorphically onto an
$m$-simplex of $X$.  Recall that the volume of $f$ is the cardinality
of $I$ (by compactness of $W$, this is finite).

The submanifold with boundary $W\setminus \bigsqcup_{i\in I} B_i$
admits a triangulation.  We apply the Cellular Approximation Theorem
\cite[Theorem 4.8]{Hatcher:book} to the restriction $f\co W\setminus
\bigsqcup_{i\in I} B_i \to X^{(m-1)}$ and obtain that it is homotopy
equivalent to a simplicial map $\bar{f} \co W\setminus
\bigsqcup_{i\in I} B_i \to X^{(m-1)}$.  Due to the homotopy
equivalence with $f$ it follows that for every $i\in I$, $\bar{f}
(\partial B_i)$ and $f (\partial B_i)$ coincide as sets.  We may then
extend $\bar{f}$ to a simplicial map $\hat{f}\co W \to
X^{(m)}$, with the same volume as $f$, homotopy equivalent to $f$, and
such that the sets $\hat{f}(B_i)$ and $f(B_i)$ coincide for every
$i\in I$.

In filling problems, when dealing with extensions of maps from
boundaries to domains, one may use an argument as above and the
version of the Cellular Approximation Theorem ensuring that if a
continuous map between CW-complexes is cellular on a
sub-complex $A$ then it is homotopic to a cellular map by a homotopy
which is stationary on $A$ \cite[Theorem 4.8]{Hatcher:book}.

For more on equivalent definitions of (filling) volumes and functions
we refer to \cite{Groft:thesis, Groft1, Groft2}.
\end{rmk}

\me

We now define a notion of filling radius.

\begin{defn}\label{defn:fillrad}
Given a hypersurface $\fh \co \calm \to X$ and a domain $\dgot \co\dd \to X$ filling it, the \emph{radius $\Rd (\dgot )$ of the domain} $\dgot$ is the minimal $R$ such that $\dgot \left(\dd^{(1)} \right)$ is in the closed tubular neighborhood of radius $R$ of $\fh \left(\calm^{(1)} \right)$.

The \emph{filling radius $\frad (\fh )$ of the hypersurface }$\fh $ is the infimum of all the filling radii of domains realizing $\fvol (\fh )$.
\end{defn}

\me

We recall two standard results that we will use later.
The first is an immediate consequence of Alexander duality and the second is the
Jordan-Schoenflies Theorem.

\begin{prop}\label{prop:jordan}
\begin{enumerate}
  \item\label{j1} Given $M$ a $k$--dimensional compact connected smooth or piecewise linear sub-manifold without boundary of $\R^{k+1}$ (or $\mathbb{S}^{k+1}$), its complement in  $\R^{k+1}$ (respectively $\mathbb{S}^{k+1}$) has two connected components.

      \me

  \item\label{j2} When $k=1$ and $M$ is a simple closed curve, there exists a homeomorphism of $\R^2$ transforming $M$ into the unit circle.
\end{enumerate}

\end{prop}

\begin{defn}\label{def:partition}
Let  $\fh \co\calm \to X$ be a $k$-dimensional $\partial V$-- hypersurface.

A \emph{partition} of $\fh$ is a finite family of $k$-dimensional
hypersurfaces $\fh_1, \fh_2,...\fh_q$, where $\fh_i\co\calm_i \to X$
are simplicial maps defined on simplicial structures
$\calm_i$ on boundaries $\partial V_i\, $, with the
following properties.  

\begin{enumerate}
  \item there exist simplicial maps $\sigma_i : V_i \to V$ which are homeomorphisms on $\Int V_i$, local homeomorphisms on $\partial V_i$ and may identify distinct simplexes of $\partial V_i$ of codimension at least $1$.  
  
 \item $V$ can be written as a (set-wise) union of $\sigma_i (V_i), i=1,2,\dots , q$; the sets $\sigma_i (\Int V_i), i=1,2,\dots , q$ are pairwise disjoint;
  \item there exists a domain filling $\fh$,  i.e. a simplicial map $\vf$ of $V$ into $X$ such that $\vf \circ \sigma_i|_{\partial V_i} =\fh_i$ for every $i\in \{ 1,2,\dots, q\}$.  
\end{enumerate}

Each of the hypersurfaces $\fh_i$ is called \emph{a contour of the
partition}.

When $q=2$, we say that $\fh$ is obtained from \emph{$\fh_1$ adjoined with  $\fh_2$}.

\end{defn}

\me

When discussing problems of filling we will often assume the existence
of a combing, as defined below.

We say that a simplicial complex $X$ has a \emph{bounded
$(L,C)$--quasi-geodesic combing}, where $L\geq 1$ and $C\geq 0$, if
for every $x\in X^{(1)}$ there exists a way to assign to every element
$y\in X^{(1)}$ an $(L,C)$--quasi-geodesic $\q_{xy}$ connecting $y$ to
$x$ in $X^{(1)}$, such that
$$
\dist (\q_{xy}(i), \q_{xa}(i))\leq L\dist (y,a)+L\, ,
$$ for all $x, y,a\in X^{(1)}$ and $i\in \R$.  Here the quasi-geodesics are assumed
to be extended to $\R$ by constant maps.

\me

The result below is well known in various contexts (CW--complexes, Riemannian geometry etc), see for instance \cite[Theorems 10.2.1, 10.3.5 and 10.3.6]{ECHLPT}. We give a sketch of proof here for the sake of completeness.

 \begin{lem}\label{lem:combing}
Let $X$ be a simplicial complex with a bounded
$(L,C)$--quasi-geodesic combing.  For every hypersurface $\fh$,
$$
\fvol (\fh ) \preceq \vol (\fh )\, \diam (\fh ).
$$
\end{lem}

\proof Consider an arbitrary
hypersurface, $\fh : \calm \to X$, where $\calm$ is a
 simplicial structure of the boundary $\partial V$. Consider the cut locus $\mathrm{Cut}$ of $V$ relative to its boundary,
and the normal map $p: V \to \mathrm{Cut}$.  Note that $V$ is
homeomorphic to the mapping cylinder of $p|_{\partial V}$ \cite[$\S 3.1.A''$]{Gromov(1983)}.

Fix a vertex $x_0$ in $\fh (\calm )$. The map $\fh$ can be extended to a quasi-Lipschitz map on $V$ (with an appropriate simplicial structure) as follows: the whole set $\mathrm{Cut}$ is sent
onto $x_0$ and for every vertex $v\in
\calm$ the geodesic $[v, p(v)]$ is sent to the quasi-geodesic in
the combing joining $x_0$ and $\fh (v)$.

The extension can be
transformed into a simplicial map as in \cite[Lemma 12]{AWP}.
Note that \cite[Lemma 12]{AWP} (unnecessarily) assumes that the complex $X$ is what they call $m$--Dehn, but do not use this hypothesis
in the proof. Nonetheless, by \cite[Theorem 1, p. 92]{AWP},
the simplicial complex $X$ that we use satisfies the $m$--Dehn
condition.

We have thus obtained a domain filling $\fh$ of
volume $\preceq\diam (\fh )\, \vol (\fh )$.
\endproof

\begin{defn}\label{def:isopFct}
The $k$\emph{--th isoperimetric function}, also known as the  $k$\emph{--th 
filling function}, of a simplicial complex $X$ is the
function $\Is_k \colon \R_+^* \to \R_+ \cup \{ \infty \}$ such that
$\Is_k (x)$ is the supremum of the filling volume $\fvol (\fh )$ over
all $k$--dimensional spheres $\fh$ of volume at most $Ax^{k}$.

The $k$\emph{--th filling radius} of the simplicial complex $X$ is the
function $\Rd_k \colon\R_+^* \to \R_+$ such that $\Rd_k (x)$ is the
supremum of all filling radii $\frad (\fh )$ over all $k$--dimensional
spheres $\fh$ of volume at most $Ax^{k}$.
\end{defn}

In what follows the constant $A>0$ from Definition
  \ref{def:isopFct} is fixed, but not made precise.  Note that
  the two filling functions corresponding to two different
  values of $A$ are equivalent in the sense of the relation 
  $\asymp$ defined in Section \ref{sect:prelim}.

\begin{rmk}
    In dealing with isoperimetric functions 
    some authors use an alternate formulation where the 
    volume is bounded above by $A r$, instead of $A r^{k}$; this yields an 
    equivalent notion (although the functions differ by a power of $k$), 
    but with the alternative definition one must modify the 
    equivalence relation to allow an additive term which is a multiple of 
    $r$ instead of $r^{k}$. 
    We use the present 
    definition because it yields a formulation consistent with the standard definition  of the divergence function that we use in this paper, see also 
    \cite{BradyFarb:div,Hindawi:div, Leuzinger:divergence, 
    Wenger:div}.
\end{rmk}

\begin{rmk}
    Note that when considering the $k$--dimensional isoperimetric 
    function (and, as we shall see below, the $k$--dimensional divergence function),  we have that  $x^{k}\asymp x$. Accordingly, in this 
    case we often represent our function by $x$, as the property of 
    having a \emph{linear filling} means the same both under our 
    choice of normalizing the volume and under the alternative 
    choice.
\end{rmk}

We may generalize the functions above, using instead of the sphere and
its filling with a ball, a hypersurface and its filling with a domain,
both modelled on a $(k+1)$--dimensional submanifold with boundary $V$
in $\R^{k+1}$.  We then define as above the filling function and
radius, denoted $\Is_V$ and  $\Rd_V$, respectively.

\begin{prop}[\cite{BBFS}]\label{prop:vsphere}
 Assume that $V$ has dimension $k+1$ at least $4$.
\begin{enumerate}
  \item Assume that $\partial V$ is connected.

  For every hypersurface $\fh \colon\partial V \to X$ there exists a simplicial map $f :\mathbb{B}^k \to \partial V$ defined on the $k$-dimensional unit ball, whose image $\mathcal{B}$ contains all the chambers that contribute to the volume of $\fh$, and a ball $\bgot :\mathbb{B}^k \to X^{(k-1)}$ filling the sphere $\fh \circ f|_{\mathbb{S}^{k-1}}$. 
  
  Therefore $\sgot$ and $\fh'$ compose a partition of $\fh$, where $\sgot $ is the sphere composed by the ball $\fh \circ f$ and the ball $\bgot $, and $\fh' $ is the hypersurface with image in $X^{(k-1)}$ obtained by adding $\fh|_{\partial V \setminus \mathcal{B}}$ to the ball $\bgot :\mathbb{B}^k \to X^{(k-1)}$, $\fh'$ with filling volume zero. In particular
\begin{equation}\label{eq:vsphere}
\fvol (\fh ) \leq \fvol (\sgot ) \, .
\end{equation}

\me

  \item If either $\partial V$ is connected or $\Is_k(x)$ is
  super-additive then the following inequality holds:
  $$
\Is_V(x) \leq \Is_k(x)\, .
$$
\end{enumerate}
\end{prop}

\proof The proof is identical to the proof in \cite[Remark 2.6, (4)]{BBFS}. \endproof

Using the terminology in the end of Definition \ref{def:partition}, we can express the result above by stating that \emph{every hypersurface $\fh$ of dimension $k\geq 3$ is obtained by adjoining to a $k$-dimensional sphere a hypersurface with image in $X^{(k-1)}$.}

\begin{rmk}
A stronger result than Proposition \ref{prop:vsphere} is Corollary 1
in \cite{Groft2} which removes the hypothesis that $\partial V$ be connected.
We do not need that generality to obtain the results of this paper.
  \end{rmk}

\medskip

By Proposition \ref{prop:vsphere}, we may assume without loss of generality that the hypersurfaces of dimension $k\neq 2$ are defined on simplicial structures of the $k$--sphere, while for $k=2$ they are defined on a simplicial structure of a surface.

\begin{lem}\label{lem:interiorvertex}
Given a hypersurface $\fh \co \calm \to X^{(m-1)}$, consider  $\frak
d\co\dd \to X^{(m)}$ a domain filling $\fh$.

Every domain ${\mathfrak{d}}'\co \dd' \to X^{(m)}$ filling $\fh$, with the same volume and image as $\frak d$ and a minimal number of chambers in $\dd'$ has the property that every chamber of $\dd'$ is either non-collapsed or it has a vertex in the boundary, $\calm^{(0)}$.
\end{lem}

\proof Consider ${\frak d}'\co\dd' \to X^{(m)}$ filling $\fh$, with the same volume and image as $\frak d$ and a minimal number of chambers in $\dd$.

Let $\{ v,w\}$ be endpoints of an edge that are not in $\calm^{(0)}$. If ${\frak d}'$ sends both $v$ and $w$ to one vertex, then by retracting the link of $\{ v,w\}$ to the link of $v$ one obtains a new domain filling $\fh$, with the same volume and image as $\frak d$ and a smaller number of chambers in its domain than ${\frak d}'$.

We may therefore deduce that every chamber in $\dd'$ is either non-collapsed or it has a vertex in $\calm^{(0)}$.\endproof

\begin{thm}[Theorem 1 and Corollary 3 in \cite{AWP}]\label{prop:qi}
 Let $X_1$ and $X_2$ be two $n$--connected locally finite CW--complexes such that for each $i\in \{1,2\}$ a group $G_i$ acts on $X_i$ cellularly and such that $X^{(n+1)}_i/G_i$ has finitely many cells.
 
 Let $V$ be an arbitrary $(n+1)$--dimensional
connected compact sub-manifold of $\R^{n+1}\, $, where $n\geq 1$, $V$ smooth or piecewise linear, and with connected boundary.

 If $X_1$ is quasi-isometric to $X_2$ then
 $$
 \Is_V^{X_1} \asymp \Is_V^{X_2}\, .
 $$
 
 In particular, for every $1\leq k\leq n$
 $$
 \Is_k^{X_1} \asymp \Is_k^{X_2}\, .
 $$
\end{thm}

    Theorem \ref{prop:qi} allows us to define $n$--dimensional filling
    functions for groups of type $\mathcal{F}_{n+1}$, up to the equivalence  relation $\asymp$.

\begin{defn}\label{def:isopFctG}
Let $G$ be a group acting properly discontinuously by simplicial
isomorphisms on an $n$-connected simplicial complex such that
$X^{(n+1)}/G$ has finitely many cells.  For every $1\leq k\leq n$, the
$k$\emph{--th isoperimetric function} of $G$ is defined to be the $k$--th
isoperimetric function of $X$. 

Likewise, the \emph{isoperimetric function of $G$ modelled on $V$} is defined to be $\Is_V^X$. 
\end{defn}

According to Theorem \ref{prop:qi} any pair of choices of
simplicial complexes as in Definition \ref{def:isopFctG} yield
filling functions which are $\asymp$-equivalent, thus the definition is
consistent.  This definition is also equivalent to the definitions appearing in
\cite{AWP,Bridson:PolynDehn,TRiley:higherIso,BBFS}.

\subsection{Filling radius estimates provided by filling functions.}

We begin by defining a simplicial object which, in all arguments on simplicial complexes, is meant to replace the intersection of a submanifold with a ball in a Riemannian manifold. 
 
Consider a simplicial map $\cgot \co \calc \to X$
representing either a domain $\dgot \co \dd \to X$ modelled on a
sub-manifold $V$ of $\R^{k+1}$, or a hypersurface $\fh \co\calm \to X$
modelled on $\partial V$.  Let $v$ be a
vertex of $\calc$ and let $r>0$.

\begin{notation}\label{not:restr1}
We denote by $\calc(v, r)$ the maximal
sub-complex of $\calc$ composed of chambers that can be connected to
$v$ by a gallery whose $1$-skeleton is entirely contained in $\cgot\iv \left (
\bar{B}(\cgot(v),r) \right)$.  Here $\bar{B}(\cgot(v),r)$ represents the
closed ball centered in $\cgot (v)$ with respect to the distance
$\dist$ on $X^{(1)}$.  Let $\partial \calc(v, r)$ denote the boundary of this
subcomplex.
\end{notation}

When $\cgot$ is a sphere or a ball, modulo some slight modifications preserving the volume, its restriction to $\calc (v,r)$
is either a domain or the whole sphere, while its restriction to $\partial \calc (v ,r)$ is a hypersurface. The same is true when $\cgot$ is a domain and $r$ is strictly less than the distance in $X^{(1)}$ between $\cgot (v)$ and $\cgot\left( \partial \calc^{(1)} \right)$.

Indeed, in case $\calc (v,r )$ is not a whole sphere, it is a subset of $\R^{k+1}$ or of $\R^k$, and small neighborhoods of points in it are either open sets in $\R^{k+1}$, respectively $\R^k$, or are homeomorphic to half-balls, or, for points in the interiors of some simplices of dimension $<k$ on the boundary, may have some other structures. It suffices to cut along all these latter simplices of dimension $<k$ (this can be done while remaining inside $\R^{k+1}$, respectively $\R^k$, and the boundary can be made into a smooth hypersurface) to obtain a complex $\calc (v,r)^{\mathrm{cut}}$ which is modelled on a smooth compact sub-manifold with boundary of $\R^{k+1}$, respectively $\R^k$. If $\g_r^{\mathrm{cut}} :\calc (v,r)^{\mathrm{cut}} \to \calc
(v,r)$ is the map gluing back along the simplices of dimension $<k$ where the cutting was done, then the restriction of $\cgot$ to $\calc
(v,r)$ must be pre-composed with $\g_r^{\mathrm{cut}}$.

\begin{notation}\label{not:restr2}
We denote by $\cgot
(v,r)$ and by $\partial \cgot (v,r)$ the domain, respectively the hypersurface, defined by restricting $\cgot$ to  $\calc (v,r)$, respectively to $\partial \calc (v,r)$, and pre-composing it with $\g_r^{\mathrm{cut}}$.
\end{notation}

A relation can be established between the filling radius and
the filling function.  Note that for $k=1$ and $\alpha =2$ this
relation was first proved in \cite[Proposition pg. 799]{Papasoglu:quadratic}.

\begin{prop}\label{prop:radius}
Let $k\geq 1$ be an integer.

If $k\neq 2$ then assume that $\Is_k (x)\leq B x^{\alpha}$ for $\alpha
\geq k$ and some constant $B>0$; while if $k=2$ then assume that for
every compact closed surface $\partial V$ bounding a 3--dimensional handlebody $V$ in $\R^3$, $\Is_V (x)\leq B
x^{\alpha}$, where $\alpha\geq 2$ and $B>0$ are independent of $V$.

Consider an arbitrary connected hypersurface $\fh\colon\calm \to X$ of dimension $k$ such that $\fvol (\fh )\geq 1$. For every filling domain $\dgot\colon\dd \to X$ realizing $\fvol (\fh )$ and such that $\dd$ has a minimal number of chambers, the following holds.
\begin{enumerate}
\item\label{rad1}  If $\alpha = k$ then $
\Rd (\dgot ) \leq C \ln \fvol (\fh )\, , $ where $C= C(A,B,k)$.

\me

\item\label{rad2}  If $\alpha > k$ then $
\Rd (\dgot ) \leq D \left[ \fvol (\fh ) \right]^{\frac{\alpha - k}{\alpha }}\, , $ where $D= D(A,B,k, \alpha )$.

\me

\item\label{rad3} If $\alpha >k$ and if for $\epsilon >0$
  small enough and for $x$ larger than
  some $x_0$, $\Is_k (x)\leq \epsilon x^{\alpha}$
  (respectively, for $k=2$, $\Is_V (x)\leq \epsilon x^{\alpha}$ and this holds for $V$ as above) then either $\Rd (\dgot ) \leq Bx_0^\alpha$ or
$$
\Rd (\dgot )\leq L \epsilon \left[ \fvol (\fh ) \right]^{\frac{\alpha - k}{\alpha }} + \dfrac{1}{\left(L \epsilon \right)^{\frac{\alpha }{\alpha -k}}} 
$$ where $L= L(A, k,\alpha)\, $.
\end{enumerate}
\end{prop}

\proof Consider an arbitrary vertex $v$ of $\dd_\vol \setminus \calm$, and a positive integer 
$$i < \dist \left( \dgot (v)\, ,\,
  \fh\left(\calm^{(1)}\right) \right)\, .
  $$  The filling domain
  $\dgot_i = \dgot (v, i )$
realizes the filling volume of the hypersurface $\fh_i =
\partial \dgot (v, i )$, by the minimality of the volume of $\dgot$.

Since $v$ is in $\dd_\vol \setminus \calm$, $\vol (\dgot_i)\geq 1$ for every $i\geq 1$.   By the isoperimetric
inequality this implies that $\vol (\fh_i)\geq 1$.

\me

When $k=1$, by Proposition
\ref{prop:jordan}\eqref{j1}, we have that $\dd(v,i)$ is homeomorphic to a disk with
holes having pairwise disjoint interiors. In particular, $\dd(v,i)$ is contained in a
simplicial disk $\dd_i\subset \dd$ whose boundary $\cals_i$ is inside
a connected component of $\partial \dd(v,i)$.  It follows that
$$
\vol (\dgot_i) \leq \vol \left( \dgot|_{\dd_i} \right) \leq B
\left( \frac{\vol \left( \dgot|_{\cals_i} \right)}{A}\right)^{\frac{\alpha}{k}}
\leq B \left( \frac{\vol \left( \fh_i \right)}{A}\right)^{\frac{\alpha}{k}}.
$$

Assume now that $k\geq 2$. By Proposition \ref{prop:jordan}\eqref{j2}, if $\partial \dd(v,i)$ is composed of several closed connected $k$--dimensional submanifolds $\sgot_1, \sgot_2,...\sgot_q$ then $\dd(v,i)$ is the intersection of connected components of $\dd \setminus \sgot_i$, one component for each $i\in \{1,2,...q\}\, $. For one  $i\in \{1,2,...q\}\, $ the connected component does not contain the boundary $\calm$. Let $\dd_i$ be that connected component. It is a domain modelled on a manifold $V$ such that $\partial V$ is connected and has a simplicial structure isomorphic to $\sgot_i$. Then
$$
\vol (\dgot_i) \leq \vol \left( \dgot|_{\dd_i} \right)
\leq \Is_V \left( \left(\frac{\vol \left( \dgot|_{\sgot_i}\right)}{A}\right)^{\frac{1}{k}}\right)
\leq \Is_V \left( \left(\frac{\vol \left( \fh_i \right)}{A}\right)^{\frac{1}{k}} \right).
$$

If $k\geq 3$ then $\Is_V  \left( \left(\frac{\vol \left( \fh_i \right)}{A}\right)^{\frac{1}{k}} \right) 
\leq \Is_k  \left( \left(\frac{\vol \left( \fh_i \right)}{A}\right)^{\frac{1}{k}} \right) $ 
by Proposition \ref{prop:vsphere}.

If $k=2$ then 
$\Is_V  \left( \left(\frac{\vol \left( \fh_i \right)}{A}\right)^{\frac{1}{k}} \right) \leq 
B \left( \frac{\vol \left( \fh_i \right)}{A}\right)^{\frac{\alpha}{k}}$ by 
hypothesis.

Thus, in all cases we obtained that
\begin{equation}\label{eq:zi}
\vol (\dgot_i) \leq B \left( \frac{\vol \left( \fh_i \right)}{A}\right)^{\frac{\alpha}{k}}.
\end{equation}

Since $i$ is strictly less than $\dist \left( \dgot (v)\, ,\,
  \fh\left(\calm^{(1)}\right) \right)\, $, the volume of $\dgot_{i+1}$ is at least the volume of $\dgot_i$
plus $\left(\frac{1}{k+1}\right)$--th of the volume of $\fh_{i}$.

Indeed, every codimension one face in $\fh_{i}$ is contained in two
chambers $\Delta , \Delta'$, such that $\Delta$ is in $\dd(v,i)$ and
$\Delta'$ is not. If $\Delta'$ is collapsed then $\dgot (\Delta' )=
\dgot (\Delta \cap \Delta' )$, whence $\Delta'$ is in $\dd(v,i)$ too.
This contradicts the fact that the codimension one face $\Delta \cap
\Delta' $ is in the boundary of $\dd(v,i)$.
It follows that $\Delta'$ is not collapsed, and it is in $\dd(v,i+1 )\setminus \dd(v,i)$.

Whence
\begin{equation}
\vol \left( \dgot_{i+1}\right) \geq \vol \left( \dgot_i \right) + C\, \vol \left( \dgot_{i+1}\right)^{\frac{k}{\alpha}}\geq \vol \left( \dgot_i \right) + C\, \vol \left( \dgot_i \right)^{\frac{k}{\alpha}}\, , \mbox{ where }
C=\frac{A}{(k+1) B^{\frac{k}{\alpha }} }\, .
\end{equation}

Part (\ref{rad1}). Assume that $k=\alpha$.  Then the above gives $\vol \left( \dgot_{i+1}\right)
\geq (1 + C )\, \vol \left( \dgot_{i}\right)$, hence by induction $\vol \left( \dgot_{i+1}\right) \geq (1+C)^{i}$.

\me

Part (\ref{rad2}). Assume that $\alpha >k \, $.  We prove by induction
on $i\leq \dist (v, \calm )$ that
\begin{equation}\label{eq:rad2}
    \vol \left( \dgot_{i}\right) \geq D i^{\frac{\alpha }{\alpha -k}}\, \mbox{ for $D$ small
    enough.  }
\end{equation}

The statement is obvious for $i=1$, and if we assume it for $i$ then
$$
\vol \left( \dgot_{i+1}\right) \geq D i^{\frac{\alpha }{\alpha -k}} + CD^{\frac{k}{\alpha
}}i^{\frac{k }{\alpha -k}} \, .
$$

Thus it suffices to prove that
$$
D\left[ (i+1)^{\frac{\alpha }{\alpha -k}} - i^{\frac{\alpha }{\alpha
-k}} \right]\leq C D^{\frac{k}{\alpha }}i^{\frac{k }{\alpha -k}} \, .
$$

A standard application of the Mean Value Theorem proves that the
latter holds if $D$ is small enough compared to $C$.

\me

Part (\ref{rad3}).  Assume $\alpha >k \, $ and moreover that for
every $x\geq x_0$, $\Is_k (x) \leq \epsilon x^{\alpha}$ (respectively
$\Is_V (x)\leq \epsilon x^{\alpha}$ for every surface, when $k=2$).

For $i\geq i_0=B x_0^\alpha$ we have that $x_i=\left(\frac{\vol\left(\fh_i \right)}{A} \right)^{\frac{1}{k}}$ is at least $x_0\, $.  Thus the
domain $\dgot_i$ filling $\fh_i$ and realizing the filling volume has
$$
\vol (\dgot_i) \leq \epsilon x_i^\alpha
$$

This implies that for $i\geq i_0$
\begin{equation}\label{eq:zi1}
\vol \left(\dgot_{i+1}\right) \geq \vol \left( \dgot_{i} \right) + C_\epsilon\, \vol \left( \dgot_{i} \right)^{\frac{k}{\alpha}}\, , \mbox{
where } C_\epsilon =\frac{A}{(k+1) \epsilon^{\frac{k}{\alpha }}
}\, .
\end{equation}

Let $D= \frac{\mu }{\epsilon^{\frac{k}{\alpha -k}}}\, ,$ where $\mu =
\left(\frac{A(\alpha -k)}{ \alpha 2^{\frac{k}{\alpha -k}} (k+1)}\right)^{\frac{\alpha }{\alpha -k}}\, $.

Consider $j_0$ large enough so that $\vol (\dgot_{j_0+1}) \geq D$.  We can take
$j_0$ to be the integer part of $D$.  We prove by induction that for
every $i\geq j_0 +1$
$$
\vol \left( \dgot_{i} \right) \geq D (i-j_0)^{\frac{\alpha }{\alpha -k}}\, .
$$

Assume that the statement is true for $i$.  According to
\eqref{eq:zi1},
$$
\vol \left(\dgot_{i+1}\right) \geq D (i-j_0)^{\frac{\alpha }{\alpha -k}} + C_\epsilon
D^{\frac{k}{\alpha }} (i-j_0)^{\frac{k }{\alpha -k}}\, .
$$

The right hand side of the inequality is larger than $D
(i+1-j_0)^{\frac{\alpha }{\alpha -k}}$ if
\begin{equation}\label{eq:de}
D\left[ (i+1-j_0)^{\frac{\alpha }{\alpha -k}} - (i-j_0)^{\frac{\alpha
}{\alpha -k}} \right] \leq C_\epsilon D^{\frac{k}{\alpha }}
(i-j_0)^{\frac{k }{\alpha -k}}\, .
\end{equation}

We may apply the Mean Value Theorem to bound the left hand side of
\eqref{eq:de} from above and write
$$
D\left[ (i+1-j_0)^{\frac{\alpha }{\alpha -k}} - (i-j_0)^{\frac{\alpha
}{\alpha -k}} \right]\leq D \frac{\alpha}{\alpha -k}
(i+1-j_0)^{\frac{k}{\alpha -k}}\leq D \frac{\alpha}{\alpha -k}
2^{\frac{k}{\alpha -k}}(i-j_0)^{\frac{k}{\alpha -k}}
$$

Thus the inequality \eqref{eq:de} holds true if
$$
D \frac{\alpha}{\alpha -k} 2^{\frac{k}{\alpha -k}} \leq C_\epsilon
D^{\frac{k}{\alpha }}\, .
$$

The value chosen for $D$ implies that we have equality.
\endproof

\me

Two types of filling function estimates, listed below, play an important part in the theory.

A simplicial complex $X$ is said to \emph{satisfy a cone-type
inequality for $k$}, where $k\geq 1$ is an integer, if for every
$k$--dimensional sphere $\fh \co\calm \to X$ its
filling volume satisfies the inequality:
\begin{equation}\label{eq:cone}
\fvol (\fh ) \preceq \vol (\fh) \diam (\fh)\, .
\end{equation}

In the above inequality the diameter $\diam (\fh )$ is the diameter of
$\fh \left( \calm^{(1)}\right)$ measured with
respect to the metric of the 1-skeleton $X^{(1)}$.

A simplicial complex $X$ is said to \emph{satisfy an 
isoperimetric inequality at most Euclidean (respectively Euclidean) for $k$} if  $\Is_{k,X} (x) \preceq x^{k+1}$ (respectively $\Is_{k,X} (x) \asymp x^{k+1}$).

An immediate consequence of Proposition \ref{prop:radius} is the
following.
\begin{cor}\label{cor:linrad}
Let $k\geq 1$ be an integer and assume that $X$ satisfies an 
isoperimetric inequality at most Euclidean for $k$,
if $k\neq 2 , $ or an inequality of the form $\Is_V (x)\leq B
x^{3}$ for every compact closed surface bounding a handlebody 
$V$, if $k=2$.

Then the filling radius, $\Rd (x)$, 
defined in Proposition \ref{prop:radius}
is bounded from above by an affine function of $x$.
\end{cor}

\section{Filling reduced to round unfolded spheres}\label{sec:roundunfold}


\subsection{Partition into round hypersurfaces}

Among the $k$--dimensional hypersurfaces, there is a particular type
for which the cone-type inequality \eqref{eq:cone} implies an isoperimetric inequality at most Euclidean.

\begin{defn}
A $k$--dimensional hypersurface
 $\fh$ is called
$\eta$--\emph{round} for a constant $\eta >0$ if $\diam (\fh) \leq \eta
\vol (\fh )^{\frac{1}{k}}\, $.
\end{defn}

Note that for $k=1$ the hypersurfaces, i.e. the closed curves, are
always round.  In what follows we extend, for $k\geq 2$, a result from
Riemannian geometry to the setting of simplicial complexes, more
precisely we prove that every sphere has a partition into round
hypersurfaces.  

\begin{prop}[partition into round hypersurfaces]\label{prop:round}
Consider an integer $k\geq 2$.  If $k\neq 3$ then assume that
$\Is_{k-1} (x)\leq B x^{k}$ for some constant $B>0$.  If $k=3$ then
assume that for every handlebody $V$ in $\R^3$, $\Is_V
(x)\leq B x^{3}$, where $B>0$ is independent of $V$.  

Then for every $\epsilon >0$ there exists a constant $\eta >0$ such that every $k$--dimensional sphere $\fh$ has a
partition with contours $\fh_1,...\fh_n$ and $\rgot$ such that $\fh_i$
are $\eta$--round hypersurfaces for every $i\in \{ 1,2,...n\}$, $\rgot
$ is a disjoint union of hypersurfaces obtained from $k$--dimensional spheres adjoined with hypersurfaces of volume and filling volume zero, and
\begin{enumerate}
  \item\label{round1} $\sum_{i=1}^n \vol (\fh_i ) \leq 2 \vol (\fh )\, $;

  \me

  \item\label{round2} $\vol (\rgot ) \leq \theta \vol (\fh ) \, ,$ where $\theta = 1-\frac{1}{6^{k+1}}$.

  \me

  \item\label{round3} $\fh_1,...\fh_n$ and $\rgot$ are entirely contained in the neighborhood $\nn_R(\fh )$, with $R= \varepsilon \vol (\fh )^{1/k}$.  
\end{enumerate}
\end{prop}

\begin{figure}[h]
\includegraphics[width=0.35\textwidth]{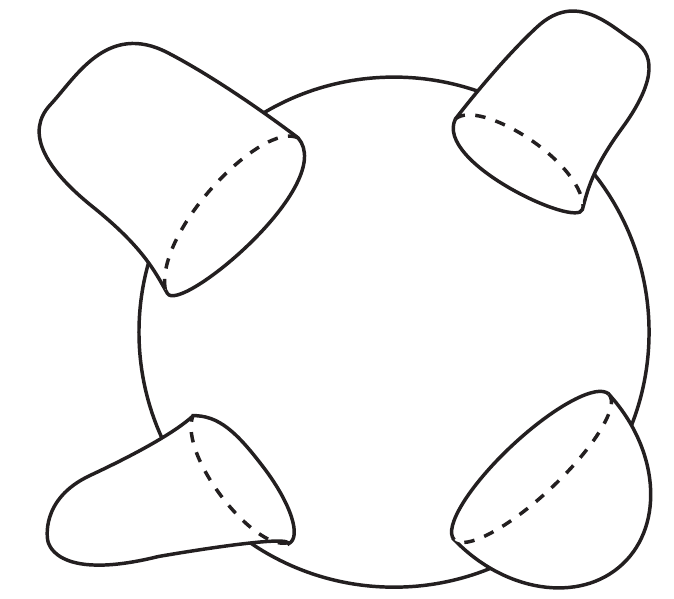}
\put(-55,50){$r$}
\caption{A sphere partitioned into round spheres and a sphere, $r$, 
of smaller area.}\label{fig:sphere_partition}
\end{figure}

\begin{rmk}\label{rem:round}
In Proposition \ref{prop:round}, for $k=2$ all the hypersurfaces $\fh_i$ are spheres, and $\rgot$ is a disjoint union of spheres. For $k\geq 3$, following Proposition \ref{prop:vsphere} one can only say that each $\fh_i$ is obtained from a $k$--dimensional sphere by adjoining it with a hypersurface of volume and of filling volume zero, and the same for all the components of $\rgot$.  
\end{rmk}

\proof Let $\fh \co \calm \to X$ be a $k$-dimensional sphere.
 We denote by $\lam >0$ a fixed small constant to be determined during the
 argument.

If the volume of $\fh$ is zero then we simply take $\rgot =\fh$. In what follows we therefore assume that $\vol (\fh )\geq 1$. For an arbitrary vertex $y\in \calm_\vol$ we define $r_0(y)$ to be the maximum of $r\geq 1$ such that $\fh (y, r)$
has volume at least $\lambda r^k$. Since $\fh (y, 1)$ contains at
least one chamber, if we assume that $\lam \leq
\frac{1}{2^{k}}\, $, we ensure that for every
$y$ the radius $r_0(y)$ is at least $2 $.

\medskip

\textsc{Step 1.} \quad Let $r_1$ be the maximum of the $r_0(y)$ for
$y\in \calm_\vol $, and let $y_1\in \calm_\vol$ be such that $r_1 =
r_0 (y_1)$.  Then consider $Y_2 = \calm_\vol \setminus \calm (y_1,
6 r_1) $, the maximum $r_2$ of the $r_0(y)$ for $y\in Y_2$,
and $y_2\in Y_2$ such that $r_2 = r_0 (y_2)$.  Assume that we have
found inductively $y_1,...,y_j$ and in $Y_{j+1} = \calm_\vol \setminus
\bigcup_{i=1}^j \calm (y_i, 6 r_i)$ consider the maximal radius
$r_0(y)$ denoted by $r_{j+1}$ and a point $y_{j+1}\in Y_{j+1}$ such
that $r_0 (y_{j+1}) = r_{j+1}$.

We thus find a sequence $y_1,...,y_N$ of vertices and a non-increasing
sequence $r_1\geq r_2\geq ..\geq r_N$ of radii, and we clearly have
that for $i\neq j$ the sets $\calm (y_i, 2 r_i)$ and $\calm (y_j,
2 r_j)$ do not contain a common chamber.  For $N$ large enough
we have that $\calm_\vol \setminus \bigcup_{i=1}^N \calm(y_i, 6 r_i)$
is empty. For each $i$, either $\vol \left( \fh (y_i, 6  r_i) \right)$ equals $\vol \left( \fh (y_i,
r_i ) \right)$ or it is strictly larger than $\vol \left( \fh (y_i,
r_i ) \right)$. In the latter case, we can write:
$$
\vol
\left( \fh (y_i, 6  r_i) \right) \leq \lambda 6^k
 r_i^k \leq 6^k \vol \left( \fh (y_i,
r_i ) \right)\, .
$$

In both cases we can write:
$$
 \vol (\fh ) \leq \sum_{i=1}^N \vol
\left( \fh (y_i, 6  r_i) \right)\leq 6^k \sum_{i=1}^N \vol \left( \fh (y_i,
r_i ) \right)\, .
$$

We may therefore conclude that the domains $\{\fh (y_i, r_i)\}_{1\leq i\leq N}$ contain
at least $\frac{1 }{6^k}$ of the volume of $\fh $.

If for some $j$ we have that $\calm_\vol \subset \calm (y_j, 6 r_j)$ then $\lam r_j^k \leq \vol (\fh ) \leq \lam (6 r_j)^k$ and this may be seen as a particular case of the above, with the set of domains $\{\fh (y_i, r_i)\}_{1\leq i\leq N}$ replaced by the singleton set $\{\fh (y_j, r_j)\}\, $.

In what follows we assume that for every $i$, $\calm_\vol $ is not contained in any $\calm (y_i, 6 r_i)$.

\medskip

\textsc{Step 2.} \quad Fix $i\in \{ 1,2,\dots , N\}$ and define the
function ${\mathcal{V}}_i (r)=\vol (\fh (y_i, r))$. By the definition of $r_i$ we have that ${\mathcal{V}}_i (r_i) \geq \lam r_i^k$ while ${\mathcal{V}}_i (r_i +1) \leq \lam (r_i +1)^k$.

Assume that $\vol (\partial \fh (y_i, r_i)) >0\, .$ We may write that
$${\mathcal{V}} (r_i +1) \geq {\mathcal{V}} (r_i) + \frac{1}{k} \vol
(\partial \fh (y_i, r_i))\, ,
$$ whence, according to the Mean Value theorem,
$$
\vol (\partial \fh (y_i, r_i))\leq \lam k^2 (r_i +1)^{k-1} \leq C \lam r_i^{k-1}
$$ where $C= 2^{k-1}k^2 $.

If $\partial \fh (y_i, r_i))$ is empty, i.e. $\fh = \fh (y_i, r_i)$ then the inequality above is automatically satisfied.

If $k=2$ then $\calm (y_i, r_i )$ is either a disk with holes with
disjoint interiors or it is the whole sphere, i.e. $\partial \calm
(y_i, r_i )$ is either empty or a union of circles.  We fill the
$\fh$-image of each circle in $X$ with a disk of area quadratic in the
length of the circle, and transform $\fh (y_i, r_i )$ into a sphere
$\fh_i$ of area at most $\lam r_i^2 + C^2 \lam^2r_i^2$.

\begin{figure}[h]
\includegraphics[width=0.35\textwidth]{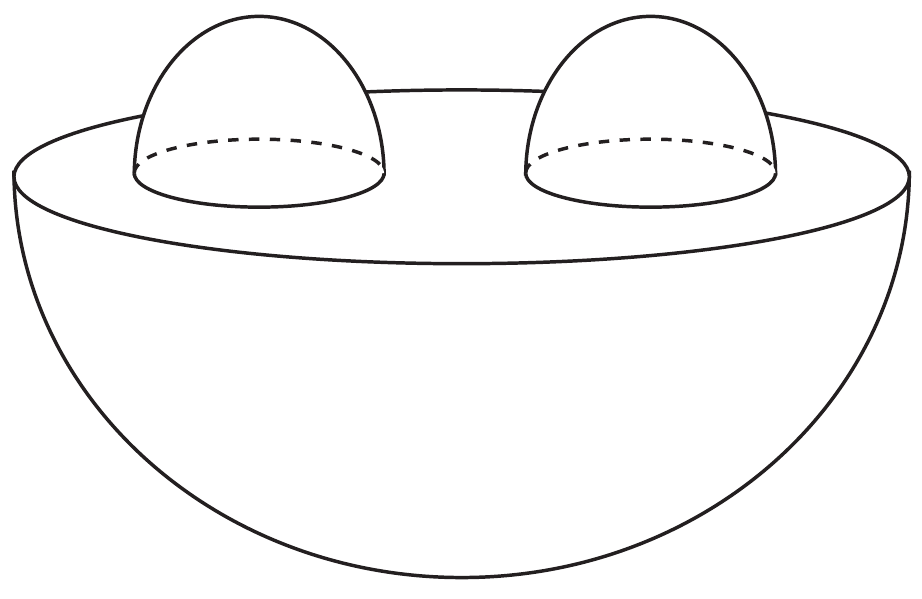}
\caption{The sphere $\fh_i$ in the case 
$k=2$.}\label{fig:sphere}
\end{figure}

The added area is at most $C^2 \lam^2r_i^2$, which for $\lam $ small enough
is at most $\frac{1}{6}\lambda r_i^2\leq \frac{1}{6}{\mathcal{V}}_i
(r_i)$. Therefore the sphere $\fh_i$ has area at most $\left(1+\frac{1}{6} \right) {\mathcal{V}}_i (r_i )$.

On the other hand, we consider the remainder of the complex $\calm \setminus \bigcup_{i=1}^N\calm (y_i, r_i )$ to which we add $\bigcup_{i=1}^N \partial \calm (y_i, r_i)$. We fill in $X$ the $\fh$-image of each circle in $\bigcup_{i=1}^N \partial \calm (y_i, r_i )$ with a disk of area quadratic in the length of the circle, and transform $\fh|_{\calm \setminus \bigcup_{i=1}^N\calm (y_i, r_i )}$ into a disjoint union of $2$-spheres, denoted by $\rgot$.

If $k\geq 3$ then $\calm (y_i, r_i )$ has boundary composed of
$(k-1)$--dimensional hypersurfaces and by using a similar argument,
the hypothesis on the $(k-1)$--filling and Proposition
\ref{prop:vsphere}(2), when $k\geq 4$, we transform each $\fh (y_i,
r_i )$ into a hypersurface $\fh_i$ of volume at most $\left(1+\frac{1}{6} \right)
{\mathcal{V}}_i (r_i )$, and the remaining complex $\calm \setminus
\bigcup_{i=1}^N\calm (y_i, r_i )$ into a disjoint union composed of a
$k$-sphere and of $k$-hypersurfaces.  Proposition \ref{prop:vsphere}(1), allows to write each of these hypersurfaces as spheres of the
same volume adjoined with $k$-hypersurfaces contained in $X^{(k-1)}$
with volume and filling volume zero.  We again denote this union by $\rgot$.

The hypersurface $\fh_i$ has volume $\asymp r_i^k$ and diameter
$\preceq r_i$, since $\fh (y_i, r_i )$ has diameter at most $2r_i$,
and for the filling domain of each component of $\partial \fh (y_i,
r_i )$ the radius is $\preceq r_i$ by Proposition \ref{prop:radius}.

In both cases $k=2$ and $k\geq 3$ we obtain that
$$
\sum_{i=1}^n \vol (\fh_i ) \leq \left(1+\frac{1}{6} \right) \sum_{i=1}^n \vol
(\fh(y_i, r_i) )\leq 2 \sum_{i=1}^n \vol
(\fh(y_i, r_i) )\, .
$$

Since for every $i\neq j$ the sets of chambers in $\fh(y_i, r_i) \cap \calm_\vol $ and respectively in  $\fh(y_j, r_j) \cap \calm_\vol $ are disjoint,
$$
\sum_{i=1}^n \vol (\fh(y_i, r_i) ) \leq \vol (\fh )\, .
$$

The union $\rgot$ is obtained by
replacing the domains $\fh (y_i, r_i )$ with unions of filling domains of their boundary components. It therefore has volume at most $\vol (\fh ) - \sum_{i=1}^N \vol \left( \fh (y_i , r_i ) \right) + \frac{1}{6} \sum_{i=1}^N \vol \left( \fh (y_i , r_i ) \right)$. We combine this with the fact that $\sum_{i=1}^N \vol \fh (y_i , r_i ) \geq \frac{1}{6^k} \vol (\fh )$ and obtain \eqref{round2}.

To prove \eqref{round3} it suffices to note that both $\rgot$ and $\fh_i$ are obtained by adding minimal volume domains filling hypersurfaces corresponding to connected components of each $\partial {\mathcal{M}} (y_i,r_i)$. Each of these components has area at most $C\lambda r_i^{k-1}$. This, the hypothesis on the filling in dimension $k-1$ and Proposition \ref{prop:radius}, \eqref{rad2}, implies that $\fh_i$ is contained in $\nn_{R_i} (\fh )$, where $R_i\leq D' \lambda^{\frac{1}{k-1}}r_i$. On the other hand $\lambda r_i^k \leq \vol (\fh )$, whence $r_i\leq \lambda^{-\frac{1}{k}} \vol (\fh )^{\frac{1}{k}}$. It follows that $R_i \leq D' \lambda^{\frac{1}{k(k-1)}} \vol (\fh )^{\frac{1}{k}}$, and the latter bound is at most $\varepsilon \vol (\fh )^{\frac{1}{k}}$, for $\lambda $ small enough. \endproof

\begin{thm}\label{thm:rounddec0}
Consider an integer $k\geq 2$.  If $k\neq 3$ then assume that
$\Is_{k-1} (x)\leq B x^{k}$ for some constant $B>0$.  If $k=3$ then
assume that for every $3$-dimensional handlebody $V$, $\Is_V
(x)\leq B x^{3}$, where $B>0$ is independent of $V$.

Then for every $\varepsilon >0$ there exists a constant $\eta >0$ such that every $k$--dimensional sphere $\fh$ has a
partition with contours $\fh_1,...\fh_n$ that are $\eta$--round hypersurfaces, and contours $\rgot_1,\dots ,\rgot_m$ that are hypersurfaces of volume and filling volume zero such that
\begin{enumerate}

\item\label{rdec1} $\sum_{i=1}^n \vol (\fh_i ) \leq 2\cdot 6^{k+1} \vol (\fh )\, .$

\medskip

\item\label{rdec2} $\fh_1,...\fh_n$ and $\rgot_1,\dots ,\rgot_m$ are  contained in the tubular neighborhood $\nn_R (\fh )$, where $R= \varepsilon \vol (\fh )^{1/k}$.
\end{enumerate}

For $k=2$ all the hypersurfaces $\fh_1,...\fh_n$ and $\rgot_1,\dots ,\rgot_m$ are spheres.  

If moreover $X$ has a bounded quasi-geodesic combing then for every $k\geq 3$ as well, $\fh_1,...\fh_n$ are $\kappa \eta$--round spheres, where $\kappa$ depends only on the constants of the quasi-geodesic combing. 
\end{thm}

\begin{rmk}\label{rmk:rounddec0}
Theorem \ref{thm:rounddec0} is true also for $k$--dimensional hypersurfaces $\fh$. Indeed, for $k\geq 3$ the statement follows from Theorem \ref{thm:rounddec0} for spheres and Proposition \ref{prop:vsphere}. For $k=2$ the same argument works and yields a decomposition with $\fh_1,...\fh_n$ and $\rgot_1,\dots ,\rgot_m$ surfaces. The only difference is that, in the proof of Proposition \ref{prop:round}, ${\mathcal{M}}(y_i, r_i)$ are subsurfaces with boundary, and after their boundary circles are filled with disks they become closed surfaces $\fh_i$.  The remaining $\rgot$, obtained from $\fh$ after performing this operation for all $i$, is likewise a disjoint union of surfaces. 
\end{rmk}

\proof  Proposition \ref{prop:round} allows to find, for an arbitrary $\varepsilon' >0$ a constant $\eta >0$ such that for every $i\geq 1$, an arbitrary
$k$-dimensional sphere $\fh$ admits a partition with contours
$\fh_1^i,...\fh_{n_i}^i , \rgot_i$ such that
$\fh_1^i,...\fh_{n_i}^i$ are $\eta$-round, moreover:

\begin{enumerate}
  \item $ \vol (\rgot_i) \leq \theta^i \vol (\fh)\, $;

  \medskip

  \item $ \sum_{j=1}^{n_i} \vol \left( \fh_{j}^i \right)\leq 2
  \sum_{\ell=0}^{i-1} \theta^\ell \, \vol (\fh)$;

  \medskip

  \item $\fh_1^i,...\fh_{n_i}^i , \rgot_i$ are contained in $\nn_{R_i} (\fh )$, where $R_i =\varepsilon' (1+\theta + \cdots + \theta^{i-1})$,
\end{enumerate}

where $\theta = 1-\frac{1 }{6^{k+1}}\, $.

Indeed, the above can be proved by induction on $i\geq 1$, where the
conclusion of Proposition \ref{prop:round} yields the initial case
$i=1$.

Assume that we found the required partition for $i$. We apply Proposition \ref{prop:round} to the disjoint union of $k$-spheres composing $\rgot_i$, once the $k$-dimensional hypersurfaces of volume and filling volume zero are removed.

We obtain $\fk_1,\dots, \fk_m$ $\eta$--round hypersurfaces and $\rgot_{i+1}$ disjoint union of $k$-spheres and of $k$-dimensional hypersurfaces of volume and filling volume zero such that

\begin{enumerate}
  \item\label{v1} $\sum_{i=1}^m \vol (\fk_i ) \leq 2 \vol (\rgot_i ) \leq 2 \theta^i \vol (\fh )\, ;$

      \me

  \item\label{v2} $\vol (\rgot_{i+1} ) \leq \theta \vol (\rgot_i ) \leq \theta^{i+1} \vol (\fh )\, ;$

       \me

  \item\label{v3} $\fk_1,\dots, \fk_m$ and $\rgot_{i+1}$ are contained in the tubular neighborhood of $\rgot_{i}$ of radius $\varepsilon' \vol (\rgot_{i})^{1/k}$.
\end{enumerate}

We then consider the set of hypersurfaces
$$
\left\{ \fh^{i+1}_1,\dots ,\fh^{i+1}_{n_{i+1}} \right\} = \left\{ \fh^i_1,\dots ,\fh^i_{n_i}, \fk_1,\dots, \fk_m  \right\}\, .
$$

\me

For large enough $i$, $\vol (\rgot_{i} )=0$, thus we obtain the required partition. According to Remark \ref{rem:round}, for $k=2$ the hypersurfaces $\fh_i$ are spheres, and $\rgot$ is a disjoint union of spheres.

For $k\geq 3$, according to Proposition \ref{prop:vsphere} each $\fh_i$ has a partition composed of a $k$--dimensional sphere  $\sgot_i$ of dimension $k$ and a hypersurface $\fh'_i$ of volume and of filling volume zero. Moreover, this partition is obtained by filling a sphere $\sigma_i$ of dimension $k-1$ and volume zero on $\fh_i$ with a ball ${\mathfrak b}_i$ also contained in $X^{(k-1)}$. The existing bounded quasi-geodesic combing allows to construct  ${\mathfrak b}_i$ so that its diameter is at most a constant $\kappa$ (depending on the quasi-geodesic constants) times the diameter of $\sigma_i$. This means that if $\fh_i$ is $\eta$--round then $\sgot_i$ is $\kappa \eta$--round.
\endproof

\begin{cor}\label{thm:rounddec}
Consider an integer $k\geq 2$.  If $k\neq 3$ then assume that
$\Is_{k-1} (x)\leq B x^{k}$ for some constant $B>0$.  If $k=3$ then
assume that for every $3$-dimensional handlebody $V$, $\Is_V
(x)\leq B x^{3}$, where $B>0$ is independent of the handlebody.

If for every round $k$--hypersurface of volume at most $Ax^{k}$ its filling
volume is at most $B x^{\alpha}$ with $\alpha \geq k$, then $\Is_k (x)
\leq \kappa B x^{\alpha},$ where $\kappa $ is a universal constant.
\end{cor}

\begin{rmk}\label{rem:thmround}
In the above, for $k=2$ it suffices to require that for every round $k$--\emph{sphere} of volume at most $Ax^{k}$ the filling
volume is at most $B x^{\alpha}$.
\end{rmk}

\proof We consider the set of contours of a partition as provided by Theorem \ref{thm:rounddec0}, 
$$
\left\{ \fh_1,\dots ,\fh_{n}, \rgot_1,\dots, \rgot_m \right\}\, .
$$

The problem of filling $\fh$ is thus reduced to the
problem of filling the round hypersurfaces $\fh_i$.  These are
filled by a volume $\leq B \sum_{j=1}^{n} Vol
(\fh_i)^{\frac{\alpha}{k} }\leq B \kappa Vol (\fh )^{\frac{\alpha}{k}}\, $.\endproof

This allows to prove a Federer-Fleming type inequality for simplicial complexes, hence for groups of finite type.

 \begin{thm}\label{thm:FedererFlemingGroups}
Assume that the simplicial complex $X$ has a bounded quasi-geodesic combing.  
\begin{enumerate}
\item\label{crd1} {\rm(Federer-Fleming inequality for
 groups).} For every $k\geq 1$, $\Is_k (x)
\preceq x^{k+1}$.  Moreover for $k=2$ the supremum of
$\Is_V (x)$ over all handlebodies $V$ is $\preceq x^{3}$.

\medskip

\item\label{crd2} Assume that for some $k\geq 2$ it is known that every round $k$--sphere of volume at most $Ax^{k}$ has filling
volume at most $B x^{\alpha}$ with $\alpha \in [k, k+1)$. Then $\Is_k (x)
\leq \kappa B x^{\alpha},$ where $\kappa $ is a universal constant.
\end{enumerate}
\end{thm}

\proof \eqref{crd1} For $k=1$ the cone filling inequality implies the quadratic 
filling inequality. For $k=2$ consider an arbitrary surface $\fh$ modelled on $\partial V$ in $\R^3$. It can be cut into round surfaces as in Proposition \ref{prop:round}. Indeed, with the notations in the proof of that Proposition, the $\fh (y_i , r_i )$ are modelled on sub-surfaces with boundary of $\partial V$. By filling the circles composing the boundary of each  $\fh (y_i , r_i )$, the surface $\fh$ is cut into round surfaces $\fh_1,...\fh_n$, and a disjoint union of surfaces $\rgot$ with $\vol (\rgot ) \leq \theta \vol (\fh )$. By iterating this decomposition, as in Theorem \ref{thm:rounddec0}, we reduce the problem of filling $\fh$ to that of filling round surfaces $\fh_1,...\fh_q$ with $\sum_{j=1}^{q} \vol (\fh_j)\leq 2 \sum_{\ell=0}^\infty \theta^\ell \, \vol (\fh)$.  Lemma
  \ref{lem:combing} allows us to deduce that $\fvol (\fh_i ) \leq B(\vol
  (\fh_i))^{\frac{3}{2}}, $ where $B$ depends only on the constants of
  the combing and on $\eta$; whence $\fvol (\fh ) \leq B'(\vol (\fh
  ))^{\frac{3}{2}}, $ where $B'$ depends only on $B$ and on $\lambda$.

By inducting on $k$, Theorem \ref{thm:rounddec} and Lemma
\ref{lem:combing} allow us to deduce the result for all $k\geq 3$.

\medskip

 \eqref{crd2} follows from \eqref{crd1} and Theorem \ref{thm:rounddec0}.  
\endproof

\subsection{Filling estimates deduced from those on round unfolded spheres}

\begin{defn}\label{def:folded}
Given a $k$--dimensional hypersurface $\fh :\calm \to X$, its $\varepsilon$--\emph{folded part at scale} $\rho$ is the set
$$
\thin_\varepsilon (\fh , \rho )= \left\{ v\in \calm_\vol \mid \exists\; r \in [1,\rho ]\mbox{ such that }  \vol (\fh (v , r)) \leq \frac{1}{2\cdot 12^{k }}\varepsilon r^k \right\} \, .
$$

A hypersurface $\fh$ with $\thin_\varepsilon (\fh , \rho )$ empty is called $\varepsilon$--\emph{unfolded at scale} $\rho$. Whenever the parameters are irrelevant, we shall simply say that a hypersurface is \emph{unfolded}.
\end{defn}

Recall that $\calm_\vol$ denotes the set of $\fh$--non-collapsed chambers, and that $v\in \calm_\vol$ means that $v$ is a vertex of one of these chambers.

\begin{prop}[removal of folded parts]\label{prop:thin0}
Consider an integer $k\geq 2$.  If $k\neq 3$ then assume that
$\Is_{k-1} (x)\leq B x^{k}$ for some constant $B>0$.  If $k=3$ then
assume that for every handlebody $V$ in $\R^3$, $\Is_V
(x)\leq B x^{3}$, where $B>0$ is independent of the handlebody.

For every $\varepsilon \in (0,1)$ small enough and every $\rho >0$ the following holds.

Assume that $k\geq 3$. Then every $k$--dimensional sphere $\fh$ has a partition with
contours $\fh_1,...,\fh_q$ and $\rgot$, where $\fh_i$ are hypersurfaces and $\rgot$ is a disjoint union of $k$--dimensional spheres adjoined with hypersurfaces of volume and filling volume zero, such that
\begin{enumerate}
  \item\label{thin1} $0< \vol (\fh_i ) \leq 2 \varepsilon
  \frac{\rho^k}{6^k }\, $;

  \me

  \item\label{thin2} $\diam (\fh_i ) \leq \frac{\sigma
  }{\varepsilon^{\frac{1}{k}}} \vol (\fh_i )^{\frac{1}{k }},$ where
  $\sigma $ depends on $k$ and the filling constants $A,B$;

  \me

  \item\label{thin3} $\sum_{i=1}^q \vol (\fh_i ) \geq \frac{1 }{
  2\cdot 12^{k}} \, \mathrm{card}\,  \thin_\varepsilon (\fh , \rho )\,
  $;

  \me

  \item\label{thin4} $\vol (\rgot ) + \frac{1}{2} \sum_{i=1}^q
  \vol (\fh_i ) \leq \vol (\fh )\, $.
  \end{enumerate}
  
Assume $k=2$. Then all the above is true for $\fh_1,...,\fh_q$ spheres and $\rgot$ a disjoint union of spheres.

Moreover $\fh$ can be taken to be a surface, in which case the above is true for $\rgot$ a disjoint union of surfaces.
\end{prop}

\begin{figure}[h]
\includegraphics[width=0.35\textwidth]{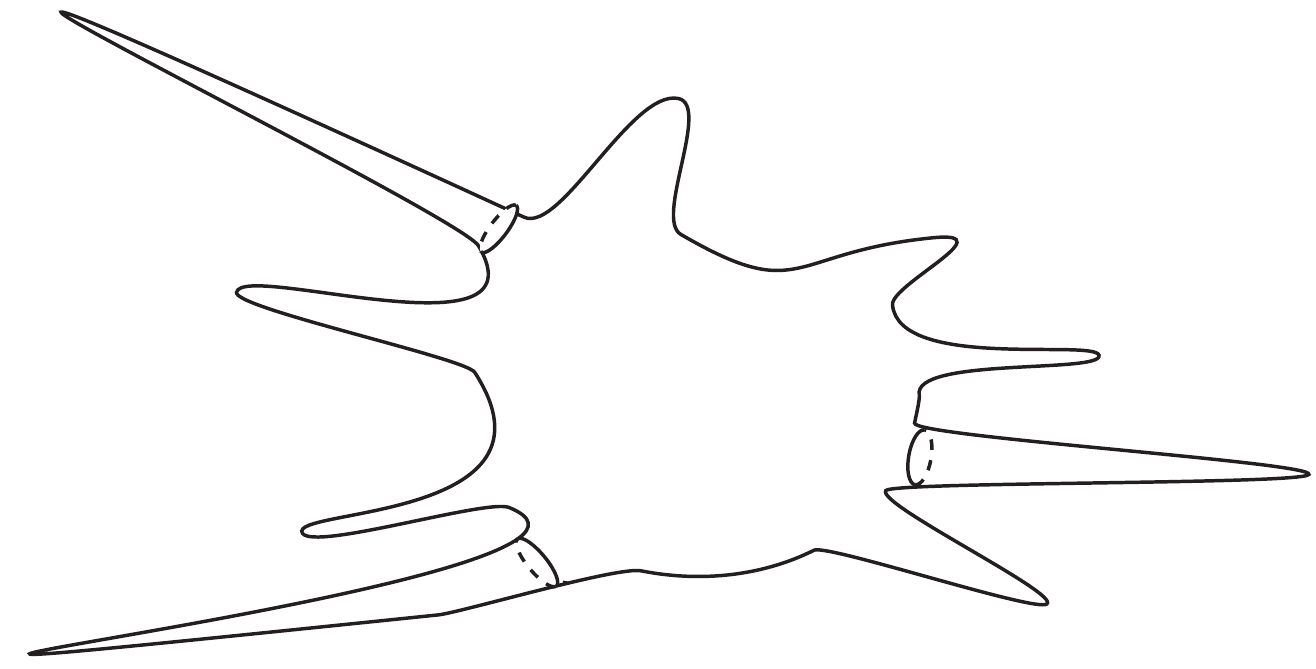}
\put(-65,25){$r$}
 \put(-120,50){$h_{1}$}
 \put(-120,-10){$h_{2}$}
 \put(0,20){$h_{3}$}

\caption{A sphere partitioned into folded spheres $h_{1}, h_{2},$ and 
$h_{3}$ and a sphere, $r$, with a smaller folded set.}
\label{fig:Sphere_folded}
\end{figure}
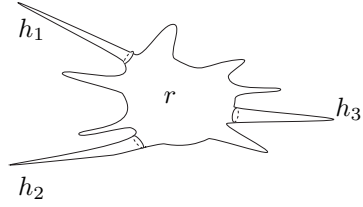

\proof We assume that $\thin_\varepsilon (\fh ,\rho )\neq \emptyset$,
otherwise we take $\rgot =\fh$ and no $\fh_1,...,\fh_q\, $.  We use
the notation introduced in \ref{not:restr1} and \ref{not:restr2}.
For every $y\in \thin_\varepsilon (\fh , \rho )$ consider $\V_y (r)=
\vol (\fh (y,r))$. Consider
$$
R_* (y) =\inf \left\{ r\in [1, \rho ]\mid \V_y (r) < \frac{1}{2
\cdot 12^{k }} \varepsilon r^k \right\} \, ,
$$ and
$$
r_*(y)= \sup \left\{ r\in [1, R_*(y) ]\mid \V_y (r) > \varepsilon r^k
\right\}\, .
$$

We assume that $\varepsilon < 1$, thus for $r=1$ we know that
$$
\vol (\fh (y , r)) > \varepsilon r^k.
$$

Note that $\varepsilon r_*(y)^k < \frac{1}{2
\cdot 12^{k }} \varepsilon R_* (y)^k$, whence $r_*(y) < \frac{R_* (y)}{12}\, $.

\begin{lem}\label{lem:ry}
 The radius $r(y) = r_*(y)+1$ satisfies the following:

\begin{enumerate}
  \item\label{b1} $\V_y (6 r(y)) \leq 12^{k} \V_y
  (r(y))$;

  \me

  \item\label{b2} $\V_y (r(y)) \leq \varepsilon r(y)^k$;

  \me

  \item\label{b3} $\vol \left( \partial \fh (y, r(y))\right)\leq C_k
  \varepsilon^{\frac{1}{k}} \V_y (r(y))^{\frac{k-1}{k}}$, where $C_k= k(k+1) 3^{k-1}$.
\end{enumerate}
\end{lem}

\proof  In what follows, for simplicity we write $r_*, R_*$ instead of $r_*(y), R_*(y)$.

The inequality
$$
\V_y (12 r_*) > 12^k \V_y ( r_*)\, ,
$$ would contradict the maximality of $r_*$. Therefore we can write
$$
\V_y (6 r(y))\leq \V_y (12 r_*) \leq 12^k \V_y ( r_*) \leq 12^k \V_y ( r(y))\, . 
$$

Property \eqref{b2} follows from the maximality of $r_*$.

We prove \eqref{b3}. We argue by contradiction and assume the inequality opposite to the one in \eqref{b3}. We can write that
$$
\V_y (r(y)+1)\geq \V_y(r(y)) +\frac{1}{k+1}\vol \left(\partial \fh (y , r ) \right) > \varepsilon r_*^k + k\varepsilon^{\frac{1}{k}} 3^{k-1}  \varepsilon^{\frac{k-1}{k}}r_*^{k-1}
$$

The right-hand side equals $\varepsilon r_*^k + k\varepsilon (3r_*)^{k-1}$, and the latter is larger than $\varepsilon (r(y)+1 )^k$, by a standard application of the Mean Value Theorem, combined with the fact that $r_* \geq 1$. This contradicts the maximality of $r_*$, hence property \eqref{b3} is true as well.\endproof

\noindent \textit{Proof of Proposition \ref{prop:thin0} continued.} \quad Consider $r_1$ to be the maximum of all $r(y)$ with $y\in \thin_\varepsilon (\fh , \rho )$ and $y_1\in \thin_\varepsilon (\fh, \rho )$ such that
$r(y_1) =r_1$.  Then consider $r_2$ to be the maximum of all the
$r(y)$ with $y\in \thin_\varepsilon (\fh, \rho ) \setminus \calm (y_1,
6 r_1)$ and $y_2$ a point in the previous set such that
$r(y_2)=r_2$.  Inductively, we find vertices $y_1,...,y_q$ and radii
$r_1,...,r_q$ and define $r_{q+1}$ as the maximum of all the $r(y)$
with $y\in \thin_\varepsilon (\fh, \rho ) \setminus \bigcup_{i=1}^q
\calm (y_i, 6 r_i)$ and $y_{q+1}$ as a point such that $r(y_{q+1})
= r_{q+1}$.

For large enough $q$ the two sequences thus constructed have the
following list of properties:
\begin{enumerate}
  \item\label{s1} $\frac{1}{2\cdot 12^{k }} \varepsilon r_i^k\leq
  \V_{y_i} (r_i) \leq \varepsilon r_i^k$;

  \me

  \item\label{s2} $r_i\leq \frac{\rho }{6}$;

  \me

  \item\label{s3} $\calm (y_i, 2 r_i)$ and $\calm (y_j, 2
  r_j)$ have no chamber in common when $i\neq j$;

  \me

  \item\label{s4} $\vol \left( \partial \fh (y_i, r_i )\right)\leq C_k \varepsilon^{\frac{1}{k}} \V_{y_i} (r_i)^{\frac{k-1}{k}}$, where $C_k= k(k+1) 3^{k-1}$;

      \me

  \item\label{s5} $\sum_{i=1}^q \V_{y_i} (r_i) \geq
  \frac{1}{2 \cdot 12^{k}} \, \mathrm{card }\,  \thin_\varepsilon (\fh, \rho )$.
\end{enumerate}

Indeed, properties (\ref{s1}), (\ref{s2}) and (\ref{s4}) follow from
the properties of $r(y)$ while (\ref{s3}) follows from the
construction of the sequences $(y_i)$ and $(r_i)$.  Property
(\ref{s5}) follows for $q$ large enough because the process can
continue until $\thin_\varepsilon (\fh, \rho ) \setminus \bigcup_{i=1}^q
\calm (y_i, 6 r_i)$ is empty, in which case
$$
{\mathrm{card}}\,  \thin_\varepsilon (\fh, \rho ) \leq 2 \sum_{i=1}^q
\V_{y_i} (6 r_i)\leq 2 \cdot  12^{k}
\sum_{i=1}^q \V_{y_i} (r_i)\, .
$$

Let $i\in \{ 1,2,\dots ,q\}$. Assume $k=2$. If $\fh$ is a sphere then $\fh (y_i , r_i)$ is a disk with other open disks removed from its interior. By filling each boundary circle
 with a quadratic area one obtains a sphere $\fh_i $ of area $\leq \V_{y_i} (r_i)(1+C_2^2 \varepsilon) $, and the remainder $\rgot$ is a disjoint union of spheres.

If we take $\fh$ to be a surface, then $\fh (y_i ,
r_i)$ is a surface with boundary. In that case, again fill each of the boundary circles
 with a quadratic area. This will yield a surface $\fh_i $ of area $\leq \V_{y_i} (r_i)(1+C_2^2 \varepsilon) $, and the remainder $\rgot$ is a disjoint union of surfaces.

 If $k=3$ then $\fh (y_i , r_i)$ is a handlebody with other
 handlebodies removed from its interior.  We fill each of the surfaces
 that compose $\partial \fh (y_i , r_i)$ with cubic volumes and
 obtain a $3$-dimensional hypersurface $\fh_i $ of volume
 $\leq \V_{y_i} (r_i)(1+C_3^{3/2} \varepsilon^{1/2})$; the remainder is
 a disjoint union of a $3$-sphere with $3$-dimensional hypersurfaces,
 and we apply Proposition \ref{prop:vsphere} (1).  A similar argument
 works in the case $k\geq 4$.

Note that for $\varepsilon$ small enough we obtain in all cases that
$$
\vol (\fh_i ) < 2 \vol \left( \fh (y_i,r_i')\right)\, .
$$

All the properties in the list are satisfied due to the properties
listed above of $y_i, r_i$ and to Proposition \ref{prop:radius}
providing the appropriate upper bound for the radius of each filling domain.
\endproof

\begin{prop}[folded-unfolded decomposition]\label{prop:thickthin}
Consider an integer $k\geq 2$.  If $k\neq 3$ then assume that
$\Is_{k-1} (x)\leq B x^{k}$ for some constant $B>0$.  If $k=3$ then
assume that for every handlebody $V$ in $\R^3$, $\Is_V
(x)\leq B x^{3}$, where $B>0$ is independent of the handlebody.

For every $\varepsilon \in (0,1)$ small enough and for every $\delta >0$ the following holds. 

Assume that $k\geq 3$. Then every
$k$--dimensional sphere $\fh$ admits a partition with contours
$\fh_1,...,\fh_q$ and $\rgot$, where $\fh_i$ are hypersurfaces and
$\rgot$ is composed of a disjoint union of $k$--dimensional spheres
which are each
adjoined with a hypersurface of volume and filling
volume zero, let $\rgot' \co \calr' \to X$ and $\rgot''\co  \calr'' \to X$
denote the simplicial maps representing the union of
$k$--dimensional spheres, respectively the union of adjoined
hypersurfaces, such that:
\begin{enumerate}
  \item ($\rgot'$ is unfolded) every vertex $v$ in $\calr'_\vol$ and every $r\in \left[ 1\, ,\, 6
  \delta \vol (\rgot)^{\frac{1}{k}} \right]$ have the property that
  $$ \vol (\rgot' (v,r )) \geq \frac{1}{2\cdot 12^k} \varepsilon r^k\, ;
  $$

  \me

  \item for every $i$, $0< \vol (\fh_i ) \leq 2 \delta^k
  \varepsilon \vol (\fh)$;

  \me

  \item $\mathrm{diam} (\fh_i ) \leq
  \frac{\sigma }{\varepsilon^{\frac{1}{k}}} \vol (\fh_i)^{\frac{1}{k}}$, where $\sigma =\sigma (A,B,k)$;

  \me

  \item $\vol (\rgot ) + \frac{1 }{2 } \sum_{i=1}^q
  \vol (\fh_i ) \leq \vol (\fh )\, $.
\end{enumerate}

Assume that $k = 2$. Then all the above is true with $\fh_1,...,\fh_q$ spheres and $\rgot = \rgot'$ a disjoint union of spheres.

Moreover one can take $\fh$ to be an arbitrary surface, in which case the statement holds with $\fh_1,...,\fh_q$ surfaces and $\rgot$ a disjoint union of surfaces.
\end{prop}

\proof The proof is by recursive applications of Proposition
\ref{prop:thin0}.  We apply Proposition \ref{prop:thin0} to the
sphere $\fh$ and $\rho = 6 \delta \vol (\fh
)^{\frac{1}{k}}$.  We obtain a partition $\rgot_1, \fh_1,...,
\fh_{q_1}$, where $\rgot_1$ is a disjoint union of $k$--dimensional
spheres each adjoined with a hypersurface of volume and filling volume zero. Let $\rgot_1'$ be the union of disjoint $k$--dimensional spheres and $\rgot_1''$ denote the union of hypersurfaces of volume and filling volume zero adjoined to each of the spheres.

At step $i$ we apply Proposition \ref{prop:thin0} to the union of spheres
$\rgot_{i-1}'$ and to $\rho_{i-1} = 6 \delta \vol (\rgot_{i-1}
)^{\frac{1}{k}}$ and we obtain a partition of $\rgot_{i-1}$ with
contours $\rgot_{i}, \fh_{q_{i-1}+1},...,
\fh_{q_{i}}$.

We have that
$$
\vol(\rgot_{i}) + \frac{1 }{2 } \sum_{j=1}^{q_i}
\vol (\fh_j) \leq \vol (\fh )\, .
$$

This in particular implies that the sum must be finite, hence the
process must stabilize at some point. Hence at some step $i$ we must find that $\thin_\varepsilon \left(\rgot_i' , \rho_i \right)$ is empty.

We take $\rgot' = \rgot_{i}'$ and
we know that $\thin_\varepsilon \left(\rgot' , 6 \delta \vol
(\rgot )^{\frac{1}{k}} \right)$ is empty.\endproof

\me

Note that ``unfolded'' does not imply ``round'', because the property of being ``unfolded'' only takes into account vertices in $\calm_\vol$. For instance, in a hypersurface $\fh : \calm \to X$ that is $\varepsilon $--unfolded at scale $6 \delta \vol (\fh)^{\frac{1}{k}}$, one can take a maximal set $S$ in $\calm_\vol$ whose image by $\fh$ is $\left[12 \delta \vol
(\fh )^{\frac{1}{k}} \right]$--separated, and its cardinality must be at most a number $N=N(\varepsilon , \delta )$. Still, the $\left[12 \delta \vol
(\rgot )^{\frac{1}{k}} \right]$--neighborhood of $\fh (S)$ only covers $\fh \left(\mathcal{\calm}_\vol \right),$ not $\fh \left(\mathcal{\calm} \right)\, $, therefore no bound on the diameter of the form $(N+2) \left[24 \delta \vol
(\fh )^{\frac{1}{k}} \right]$ can be obtained.

\me

\begin{prop}[decomposition of unfolded into unfolded and round]\label{prop:thickround}
Assume that the simplicial complex $X$ has a bounded quasi-geodesic combing and let $k\geq 2$ be an integer.

For every two numbers $\varepsilon$ and $\delta$ in $(0,1)$ there exists $N=N(\varepsilon, \delta , k)$ such that the following holds. 

If $k\geq 3$ then consider an arbitrary disjoint union of $k$--dimensional spheres $\rgot :\calr \to X$ such that $\vol (\rgot ) \geq V_0$ for some large enough constant $V_0$.

 Assume that $\rgot$ is $\varepsilon $--unfolded at scale $6 \delta \vol (\rgot)^{\frac{1}{k}}$, in the sense that
\begin{center}
(*) every vertex $v$ in $\calr_\vol$ and every $r\in \left[ 0\,
,\, 6 \delta \vol (\rgot)^{\frac{1}{k}} \right]$ have the property
that
  $$ \vol (\rgot (v,r )) \geq \frac{\varepsilon}{2\cdot 12^k}  r^k.
  $$
\end{center}
Then $\rgot$ has a partition into a union $\rgot_0$ of hypersurfaces of volume zero and $m$ spheres of dimension $k$, $\rgot_i :\calr^{(i)} \to X$ with $i\in \{1,2,...,m\}$ and $m\leq N\, $, such that
\begin{enumerate}
  \item\label{vol1} for every $i\in \{1,...,m\}$, every vertex $v$ in $\calr^{(i)}_\vol$ and every $r\in \left[ 0\, ,\, 6 \delta \vol (\rgot)^{\frac{1}{k}} \right]$ have the property that
  $$ \vol (\rgot_i (v,r )) \geq \frac{1}{2\cdot 12^k} \varepsilon r^k\, ;
  $$

  \me

  \item\label{vol2} $\diam (\rgot_i ) \leq \frac{\kappa }{\varepsilon \delta^{k-1}} \vol (\rgot)^{\frac{1}{k}}\, ,$ where $\kappa $ depends only on the constants $L$ and $C$ of the combing;

  \me

  \item\label{vol3} $\frac{\varepsilon \delta^k}{2} \vol (\rgot) \leq \vol (\rgot_i) \leq \vol (\rgot)$.
\end{enumerate}

If $k=2$ then $\rgot_0$ does not appear.

Moreover, if $\rgot $ is taken to be a disjoint union of surfaces then the statement holds with $\rgot_i,\, i\in \{1,\dots ,m\}$, surfaces, and $\rgot_0$ again does not appear.  
\end{prop}

\proof If $\vol (\rgot ) =0$ then take $\rgot_0 =\rgot\, $. We assume that $\vol (\rgot ) > 0$.

Let $q$ be the maximal number of connected components of $\rgot$ with positive volume. Property (*) implies that $q$ has a uniform upper bound depending on $\varepsilon$. Thus, without loss of generality, we may assume that $\rgot$ is one $k$--dimensional sphere.

\me

Let $v$ be an arbitrary vertex in $\calr_\vol$. Let $W$ denote $24
  \delta \vol (\rgot)^{\frac{1}{k}}$.

  We choose the constant $V_0$ in the hypothesis large enough so that $W\geq 6$.

  Suppose that for every $r\in \left[W, \frac{17}{\varepsilon \delta^k} W\right]$, $\calr (v, r+W ) \setminus \calr (v, r )$ contains a chamber from $\calr_\vol$.

  We divide $\left[W, \frac{17}{\varepsilon \delta^k} W\right]$ into consecutive intervals with disjoint interiors $[r_1-2W, r_1+2W],...,[r_q-2W, r_q+2W]$. Each $\calr (v, r_i+W ) \setminus \calr (v, r_i )$ contains a chamber from $\calr_\vol \, ,$ and for a vertex $v_i$ of that chamber the chambers in $\calr (v_i, W)$ are all contained in $\calr (v, r_i +2W ) \setminus \calr (v, r_i -2W )$. Therefore the sets of chambers in $\calr (v_i, W)$ are pairwise disjoint, and the cardinality of each set intersected with $\calr_\vol$ is at least $\frac{\varepsilon}{2\cdot 12^k}W^k$ by property (*).

  It follows that $\vol (\rgot )$ is at least $\frac{\frac{17}{\varepsilon \delta^k} W-W}{4W} \frac{\varepsilon}{2\cdot 12^k}W^k\geq \frac{\frac{16}{\varepsilon \delta^k} W}{8W} \varepsilon \delta^k \vol (\rgot )\geq 2\vol (\rgot )$. This is a contradiction.

  We conclude that there exists $r\in \left[W, \frac{17}{\varepsilon
  \delta^k} W\right]$ such that $\calr (v, r+W ) \setminus \calr (v,
  r )$ does not contain a chamber from $\calr_\vol$.  The boundary of
  $\calr_1 = \calr (v, r+W/2 )$ has no volume.  Assume for a
  contradiction that there exists a $(k-1)$-chamber $H$ of the
  boundary which is sent by $\rgot$ onto a $(k-1)$-dimensional simplex of $X$. Then there
  exists a non-collapsed chamber containing $H$, and hence contained in
  $\calr (v, r+W )$, but not in $\calr (v, r+W/2 )$.  This contradicts
  the hypothesis.

 It follows that the boundary of $\calr_1 = \calr (v, r+W/2 )$ has no
  volume.  This boundary is
  composed of hypersurfaces contained in $X^{(k-2)}$, each can be
  filled with a domain inside $X^{(k-1)}$ using the combing (see Lemma
  \ref{lem:combing}).  For $k\geq 3$, $\rgot_1$ becomes a sphere
  $\rgot_1'\co \calr_1'\to X$ with a hypersurface $\rgot_1''$ of volume
  and filling volume zero adjoined to it, and $\rgot \setminus
  \rgot_1$ becomes $\overline{\rgot}_1$, a disjoint union of spheres
  with hypersurfaces of volume and filling volume zero adjoined to
  them.  We denote the disjoint union of spheres by
  $\overline{\rgot}_1'\, $ and the union of hypersurfaces by
  $\overline{\rgot}_1''\, $.

  Assume $k=2$. Assume $\rgot$ is a disjoint union of spheres. Then $\calr_1$ is topologically a disc with open discs removed from its interior. The hypothesis that the boundary of $\calr_1$ has no
  volume implies that each boundary circle is sent onto one vertex. Therefore, by replacing $\calr$ with the quotient simplicial complex in which the above mentioned boundary circles become points, $\rgot_1$ becomes a sphere, while $\rgot$ restricted to $\calr \setminus \calr_1$ becomes a disjoint union of spheres~$\overline{\rgot}_1'$.

Assume now that $\rgot $ is a disjoint union of surfaces. Then $\calr_1$ is topologically a surface with boundary. An argument as above implies that by replacing $\calr$ with the quotient simplicial complex in which each boundary circle becomes a point, $\rgot_1$ becomes a closed surface $\rgot_1'$, while $\rgot$ restricted to $\calr \setminus \calr_1$ becomes a disjoint union of closed surfaces~$\overline{\rgot}_1'$.

We prove that, for all $k\geq 2$, both $\rgot_1'$ and $\overline{\rgot}_1'$ satisfy the property (*).
 Indeed, the part added to $\rgot_1$ to become $\rgot_1'$, and to $\rgot$ restricted to $\calr \setminus \calr_1$ to become the union of $\overline{\rgot}_1'$ with $\overline{\rgot}_1''$ does not contribute to the volume, it suffices therefore to check (*) for vertices $a$ in $\calr_1^\vol$, respectively in $\left( \calr \setminus \calr_1\right)^\vol$.

  \me

  Let $a\in \left(\calr_1\right)_\vol$. Then $a$ is in a non-collapsed chamber in $\calr (v,r)$, since $\calr_1 \setminus \calr (v,r)$ does not contain  non-collapsed chambers. It follows that $\calr \left(a, \frac{W}{4} \right) \subseteq \calr \left(v, r+ \frac{W}{2} \right)$. Thus, for every $t\in \left[0, \frac{W}{4}\right]$, $\rgot_1 (a,t) =\rgot (a,t)\, ,$ hence (*) is satisfied for $\rgot_1$.

 Likewise, let $a$ be in  $\left( \calr \setminus \calr_1\right)_\vol$. Hence $a\in \calr \setminus \calr (v, r+W)$.
 If $\calr \left( a, \frac{W}{4} \right)$ would intersect $\calr_1$ in a non-collapsed chamber then $a$ would be in $\calr (v, r+W)$, a contradiction.
Thus $\calr \left( a, \frac{W}{4} \right)_\vol$ is in $\calr \setminus \calr_1$ and as before we conclude that $\rgot$ restricted to $\calr \setminus \calr_1$ satisfies (*).

  In particular the volume of $\rgot_1'$ is at least $\frac{\varepsilon \delta^k}{2} \vol (\rgot)$.

   The diameter of $\rgot_1'$ is at most $\kappa_1 (r+W),$ where $\kappa_1$ depends on the constants $L$ and $C$ of the combing. An upper bound for $r+W$ is $\left( \frac{17}{\varepsilon\delta^k}+1 \right)W \leq \frac{18}{\varepsilon\delta^k} 24 \delta \vol (\rgot )^{\frac{1}{k}}= \frac{\kappa_2}{\varepsilon\delta^{k-1}} \vol (\rgot )^{\frac{1}{k}}\, .$

   We thus obtain the estimate in \eqref{vol2} for the diameter of $\rgot_1'$.

  We repeat the argument for $\overline{\rgot}_1'\, $, and find $\rgot_2$ etc. The process must stop after finitely many steps because of inequality (\ref{vol3}).\endproof

  We now state and deduce the main result of this section:
  
 \begin{thm}\label{thm:roundThick}
Let $X$ be a simplicial complex with a bounded quasi-geodesic combing, and let $\eta , \varepsilon$ and $\delta$ be small enough positive constants.

\begin{enumerate}
\item Let $k\geq 2$ be an integer. If every $k$--dimensional sphere of volume at most $Ax^{k}$ that is
$\eta$--round and $\varepsilon$---unfolded at scale $\delta x$, in the sense of \eqref{vol1}, has filling volume at most $B x^{\alpha}$ with $\alpha \geq k$, then $\Is_k (x) \leq A x^{\alpha},$ where $A =A (\eta , \varepsilon , \delta )$.

\item If every (closed) surface of volume at most $Ax^{2}$ that is
$\eta$--round and $\varepsilon$---unfolded at scale $\delta x$ has filling volume at most $B x^{\alpha}$ with $\alpha \geq 2$ and $B$ independent of the genus, then $\Is_V (x) \leq A x^{\alpha},$ for every handlebody $V$, where $A =A (\eta , \varepsilon , \delta )$.

\end{enumerate}
 \end{thm}

 \proof Proposition \ref{prop:thickround} implies that every disjoint union of $k$--dimensional spheres (respectively, surfaces), that is of volume at most $Ax^{k}$ and $\varepsilon$---unfolded at scale $\delta x$, has filling volume at most $BN x^{\alpha}$, where $N =N (\varepsilon , \delta ,k )$.

We prove by induction on $n$ that $k$--dimensional spheres (respectively surfaces) $\fh $ of volume at most $2^n$ satisfy 
$$ \fvol (\fh )\leq A \left( \vol (\fh ) \right)^{\frac{\alpha }{k}}$$ for $A$ large enough.

Indeed, for $k=2$ it suffices to use the decomposition in Proposition \ref{prop:thickthin}, and to apply the inductive hypothesis to each $\fh_i$. For $k\geq 3$, according to Proposition \ref{prop:vsphere}, each $\fh_i$ has a partition composed of a sphere $\fh_i'$ and a hypersurface $\fh_i''$ of volume and filling volume zero. We apply the inductive hypothesis to each $\fh_i'$.  
 \endproof
 
\comment

\begin{cor}\label{cor:dimAsCone}
Let $X$ be a simplicial complex with a bounded quasi-geodesic combing,
such that for every asymptotic cone of $X$ the maximal dimension of
locally compact subsets in it is $m$.  If $k\geq m$ then $\Is_k (x) =
o(x^{k+1})$.
\end{cor}

\proof According to Theorem \ref{thm:roundThick}, it suffices to
prove the statement for a filling function defined only for spheres (respectively surfaces) that are round and unfolded.  We argue for a contradiction and assume
that there exists a sequence $\rgot_n$ of $k$-dimensional spheres for $k\geq 3$ (respectively of surfaces for $k=2$) that
are round and unfolded, of volume $\asymp x_n^{k}$ and of filling volume at
least $\lambda x_n^{k+1}$, where $\lambda$ is a positive constant and
$x_n \to \infty\, $.  Let $\dgot_n$ be filling $(k+1)$-dimensional
balls (respectively filling handlebodies) realizing $\fvol (\rgot_n )$ and with a minimal number of chambers in the domain.

The argument in \cite[pp. 263--264]{Wenger:asrk} with
$T_n= \rgot_n$ and $S_n = \dgot_n$  can then be applied to show
that an asymptotic cone
of $X$ contains a compact subset of dimension $k+1$, a
contradiction.\endproof

As a consequence of the decompositions described previously, we obtain
that hyperbolic groups have linear filling in every dimension also for
the filling functions that we use.  The same result has been proven
for the filling functions in terms of Lipschitz spheres and Lipschitz
balls in \cite{Lang:linearHyp}.

\begin{cor}
Let $G$ be a finitely generated Gromov-hyperbolic group.  For every
$k\geq 1$, $\Is_k^G (x) \asymp x$.
\end{cor}

\proof We can consider a Rips complex $X$ for the group $G$. It suffices to prove the statement for the filling functions defined for spheres that are round and unfolded. Using the same argument as in \cite{Lang:linearHyp}, the problem reduces to the problem of estimating the filling for such spheres in a real hyperbolic space. The linear filling is then given by the same hyperbolic cone as constructed in \cite{Lang:linearHyp}.\endproof

\endcomment

\section{Divergence.}\label{section:div}

We begin by recalling a few facts about the \emph{$0$--dimensional divergence}. There are several,
essentially equivalent versions of $0$--dimensional divergence. The
first careful study of the notion was undertaken by S.~Gersten \cite{Gersten:divergence, Gersten:divergence3}. The main reference for the first part of this section is \cite[$\S
3.1$]{DrutuMozesSapir}.

Let $X$ denote a geodesic metric space, quasi-isometric to a
one-ended
finitely generated group (this assumption can be replaced by a weaker
technical hypothesis called $(\rm{Hyp}_{\kappa ,L})$, but we do not need
that generality here). Also, we fix a constant $0<\delta<1$.

For an arbitrary triple of distinct points $a,b,c\in X$ we define
$\dv(a,b,c;\delta)$ to be the infimum of the lengths of
paths connecting $a, b$ and avoiding the ball centred at $c$ and of radius $\delta \cdot
\dist(c,\{a,b\}$.  If no such path exists, define
$\dv(a,b,c;\delta)=\infty$.

    The {\em divergence function} $\Dv^{X}(n ,\delta)$ of
    the space $X$ is defined as the supremum of all finite numbers
    $\dv(a,b,c;\delta)$ with $\dist(a,b)\le n$.  When the
    ambient space is clear we omit $X$ from the notation.

It is proven in \cite[Lemma 3.4]{DrutuMozesSapir} that, as long as
$\delta$ is sufficiently small and $n$ sufficiently large, $\Dv^{X}(n ,\delta)$ is always defined and, by construction, takes only finite values.
In \cite[Corollary 3.12]{DrutuMozesSapir}, it is proven that, up to
the equivalence relation $\asymp$ (which in this case means up to 
affine functions), the
various standard notions of $0$--dimensional divergence agree and
that the $\asymp$ class of the divergence function is invariant under
quasi-isometry. As our main results are only about $\asymp$ classes 
of functions, it is no loss of generality to assume that the 
value of $n$ used in the above function is taken sufficiently large 
so that $\Dv^{X}(n ,\delta)$ is defined; hence we will make this 
assumption for $n$ for the remainder of the paper.

An important feature of the $0$--dimensional divergence function is its
relationship to the topology of asymptotic cones as described in \cite[Proposition 1.1]{DrutuMozesSapir}. 

\comment
\begin{prop} [Proposition 1.1 in
    \cite{DrutuMozesSapir}]\label{lm02}
Let $X$ be a geodesic metric space.
\begin{itemize}
\item[(i)] If there exists
$\delta \in (0,1)$ such that 
$\Dv(n;\delta) = O(n)$ then every
asymptotic cone of $X$ is without cut-points.

    \me

\item[(ii)] Assume that for some constants $c, \lambda, \kappa$, every point in $X$
is at distance at most $c$ from a bi-infinite $(\lambda ,
\kappa)$--quasi-geodesic.

If every asymptotic cone of $X$ is without cut-points then for every
$0<\delta <\frac{1}{54}$, $\Dv(n;\delta) = O(n)$.
\end{itemize}
\end{prop}

In \cite[Proposition 1.1]{DrutuMozesSapir}, in the statement (ii) it
is assumed that $X$ is periodic.  However, the periodicity condition
is only used to ensure the existence of bi-infinite quasi-geodesics
near every point.
\endcomment

We now proceed to discuss an extension to higher dimensions of the
divergence function defined above.
In an arbitrary dimension, the divergence may be seen as a filling
function when moving towards infinity (e.g. when moving closer and closer to the boundary, if a boundary $\partial_\infty$ can be defined).

The notion has been used mostly in the setting of non-positive
curvature (see e.g. \cite{BradyFarb:div} for a version defined for Hadamard manifolds). Among other things, in a Hadamard space the divergence can
distinguish the rank.  Indeed, for a symmetric space $X$ of
non-compact type, $\Dv_k $ grows exponentially when $k=\Rk (X)-1$
\cite{BradyFarb:div,Leuzinger:divergence}, while when $k\geq \Rk (X)$
the divergence $\Dv_k = O(x^{k+1})$ \cite{Hindawi:div}.  More
generally, for a cocompact Hadamard space $X$ and for a homological
version of the divergence, defined in terms of integral currents, if
$k=\Rk (X)-1$ then $\Dv_k \succeq x^{k+2}\, $, while if $k\geq \Rk
(X)$ then $\Dv_k \preceq x^{k+1}$ \cite{Wenger:div}.

In what follows, we define a version of the higher dimensional
divergence functions in the setting of simplicial complexes, in
particular of groups of type $\mathcal{F}_n$.  Therefore, we fix an
$n$-connected simplicial complex $X$ of dimension $n+1$
which is the universal cover of a compact simplicial complex $K$ with
fundamental group $G$.  Recall that we assume edges in $X$ to be of
length one, and that we endow $X^{(1)}$ with the shortest path metric.
We also fix a constant $0<\delta<1$. 
Given a vertex $c$ in $X$, a $k$--dimensional hypersurface $\fh \colon \calm
\to X$ modelled on $\partial V$ such that $k\leq n-1$, and a number $r>0$ that is at most the
distance from $c$ to $\fh (\calm^{(1)})$, \emph{the 
divergence} of this quadruple, denoted $\dv (\fh, c ;r, \delta )$, is 
the infimum of all volumes of
domains modelled on $V$ filling $\fh$ and disjoint from $B(c, \delta r)$.  If
no such domain exists then we set $\dv(\fh, c ;r, \delta )=
\infty$.

\begin{defn}\label{def:kdiv}
Let $V$ be a manifold as described in Convention \ref{cvn:V}.

The \emph{divergence function modelled on $V$} of the complex $X$, denoted $\Dv_V (r,
\delta )$, is the supremum of all finite values of 
$\dv (\fh, c; r,
\delta )$, where $\fh$ is a hypersurface modelled on $\partial V$ with the distance from $c$ to $\fh (\calm^{(1)})$ at least
$r$ and $\vol (\fh )$ at most $A r^{k}$.

When $V$ is the $(k+1)$--dimensional unit ball, $\Dv_V (r,
\delta )$ is denoted $\Dv^{(k)} (r,
\delta )$, and it is called the $k$--\emph{dimensional divergence function} (or the \emph{$k$--th divergence function}) of $X$.
\end{defn}

In the above, as for the filling functions, we fix the
constant $A>0$
 once and for all, and we do not mention it anymore. 
We note that when $k=0$, the volume condition is vacuously 
satisfied and thus the above notion coincides with that of the 
$0$--dimensional divergence given previously.

\me

An immediate consequence of Proposition~\ref{prop:radius} is that when
the isoperimetric function $\Is_V$ is smaller that the Euclidean one,
the divergence function $\Dv_V$ coincides with the isoperimetric function.

\begin{prop}\label{lem:subEucl}
Let $V$ be a manifold as described in Convention \ref{cvn:V}, and let $\varepsilon$ and $\delta$ be small enough positive constants.
Assume that $\Is_V (x) \leq \varepsilon x^{k+1}$.
Then $\Dv_V (x , \delta ) = \Is_V (x)$ for every $x$ large enough.	
\end{prop}

\proof Proposition \ref{prop:radius}, \eqref{rad3}, implies that
$\Rd_V (x) \leq 2L\varepsilon x$ for $x\geq x_\varepsilon$.  It follows that for every $k$--dimensional
hypersurface $\fh \colon \calm \to X$ of area at most $A x^{k}$, there
exists a filling domain $\dgot\colon\dd \to X$ realizing $\fvol (\fh)$
and with the image $\dgot (\dd^{(1)})$ entirely contained in a tubular
neighborhood of $\fh (\calm^{(1)} )$ of radius $2L\varepsilon x$.  Therefore, if
$\fh (\calm^{(1)} )$ is disjoint from a ball $B(c,x)$
then, for $x$ large enough, $\dgot (\dd^{(1)})$ is disjoint from the
$(\delta x)$--ball around~$c$, provided that $2L\varepsilon + \delta <1$.
\endproof

More importantly, the cutting arguments that we have described previously allow to reduce the problem of estimating the divergence to hypersurfaces that are round, in the particular case when a bounded combing exists.

\begin{thm}\label{thm:DivRound}
Assume that $X$ is a simplicial complex of dimension $n$ endowed with a bounded quasi-geodesic combing. Let $V$ be a $(k+1)$--dimensional
connected compact sub-manifold of $\R^{k+1}\, $ with connected boundary, where $2\leq k\leq n-1$.

For every $\varepsilon >0$ there exists $\eta >0$ such that the following holds. 

Consider the restricted
divergence function $\Dv^r_V(x, \delta )$, obtained by taking the supremum only over hypersurfaces modelled on $\partial V$ that are
$\eta$--round, of volume at most $2A x^{k}$ and situated outside balls of radius $x$. 

Assume that $\Dv^r_V(x, \delta ) \leq Br^\beta$ for some $\beta \geq k+1$ and $B>0$ universal
constant.  Then the general divergence function $\Dv_V (x, \delta (1- \varepsilon ) )$ is at most $B'r^\beta$ for some $B'>0$ depending on $B, \varepsilon , \eta $ and $X$.
\end{thm}

\proof Let $\fh $ be a hypersurface modelled on $\partial V$, of volume of most $Ax^k$ and with image outside the ball $B(c, x)$.

According to Theorem \ref{thm:rounddec0} and Remark \ref{rmk:rounddec0}, for every $\epsilon$ there exists $\eta = \eta (\varepsilon ) $ such that $\fh $ can be decomposed into $\eta$--round spheres $\fh_1,\dots, \fh_n$ and hypersurfaces $\rgot_1,\dots , \rgot_m$ of volume and filling volume zero, all of them contained in $\nn_{\varepsilon x }(\fh )$.

All the spheres $\fh_1,\dots, \fh_n$ have area at most $2A x^{k}$. 

We put aside the spheres $\fh_i$ that have volume at
most $\varepsilon x^{k}$. For
all these, by Proposition \ref{prop:radius}, the filling radius is at
most $L\varepsilon x$ for $x$ large enough, where $L$ is a universal constant, therefore they can be filled in
the usual way outside $B(c, \delta x)$ if $\varepsilon$ is small
enough.

In the end we obtain finitely many round components $\fh_1,\dots,
\fh_N$ of volume at least $\varepsilon x^{k}$ and at most $2Ax^{k}$, where
$N=O\left( \frac{A}{\varepsilon } \right) $, by Theorem
\ref{thm:rounddec0}, \eqref{rdec1}. Moreover $\fh_1,\dots,
\fh_N$ have images outside $B(c, (1-\varepsilon )x)$  If $\fh_1,\dots, \fh_N$ can be filled with a
domain of volume at most $Bx^\beta$ outside $B(c, \delta (1-\varepsilon ) x)$, then $\fh$ can be filled outside the same ball with a
volume at most $NB x^\beta + A' x^{k+1}$, where the second term comes
from the filling of the hypersurfaces with small volume that we had put aside, and from Theorem \ref{thm:rounddec0}.  \endproof


\newcommand{\etalchar}[1]{$^{#1}$}
\def\cprime{$'$} \def\cprime{$'$} \def\cprime{$'$} \def\cprime{$'$}
  \def\cprime{$'$} \def\cprime{$'$}
\providecommand{\bysame}{\leavevmode\hbox to3em{\hrulefill}\thinspace}
\providecommand{\MR}{\relax\ifhmode\unskip\space\fi MR }
\providecommand{\MRhref}[2]{%
  \href{http://www.ams.org/mathscinet-getitem?mr=#1}{#2}
}
\providecommand{\href}[2]{#2}


\begin{thebibliography}{ABDDY}

  \bibitem[ABDY]{ABDY:dehndifferent}
  A.~Abrams, N.~Brady, P.~Dani, and R.~Young, \emph{Homological
  and homotopical Dehn functions are different},
  Proc. Natl. Acad. Sci. USA, \textbf{110} (2013), 19206--19212.

\bibitem[AK]{AmbrosioKirchheim}
L.~Ambrosio and B.~Kirchheim, \emph{{Currents in metric spaces}}, Acta Math.
  \textbf{185} (2000), 1--80.

\bibitem[AWP]{AWP}
J.M. Alonso, X.~Wang, and S.J. Pride, \emph{Higher-dimensional isoperimetric
  (or Dehn) functions of groups}, J. Group Theory \textbf{2:1} (1999), 81--122.

\bibitem[BBFS]{BBFS}
N.~Brady, M.~Bridson, M.~Forester, and K.~Shankar, \emph{Snowflake groups,
  {P}erron-{F}robenius eigenvalues and isoperimetric spectra}, Geom. Topol.
  \textbf{13} (2009), no.~1, 141--187.

\bibitem[BBI]{Buragos:course}
D.~Burago, Y.~Burago, and S.~Ivanov, \emph{A course in metric geometry},
  Graduate Studies in Math. vol.~33, American Mathematical Society,
  Providence, RI, 2001.


\bibitem[Beh]{Behrstock:asymptotic}
J.~Behrstock, \emph{{Asymptotic geometry of the mapping class group and
  {T}eichm\"{u}ller space}}, {Geom. Topol.} \textbf{10} (2006), 1523--1578.


  \bibitem[BD1]{BehrstockDrutu:thick2}
   J.~Behrstock and C. Dru\c{t}u,
   \emph{Divergence, thick groups, and short conjugators}, \textsc{arXiv:1110.5005}.

  \bibitem[BD2]{BehrstockDrutu:HigherFilling2}
   J.~Behrstock and C. Dru\c{t}u,
   \emph{Higher dimensional divergence for
mapping class groups and $CAT(0)$ groups}, \textsc{arXiv:1305.2994}.


  \bibitem[BD3]{BehrstockDrutu:Linear}
   J.~Behrstock and C. Dru\c{t}u,
   \emph{Isoperimetry above the rank for groups with bicombings},
   preprint in preparation.

  \bibitem[BEW]{BestvinaEskinWortman}
  M. Bestvina, A. Eskin and K. Wortman \emph{Filling boundaries of coarse manifolds in semisimple and
	      solvable arithmetic groups}, \textsc{arXiv:1106.0162}, 2011.


\bibitem[BF]{BradyFarb:div}
N.~Brady and B.~Farb, \emph{Filling invariants at infinity for manifolds of
  non-positive curvature}, Trans. Amer. Math. Soc. \textbf{350} (1998), no.~8,
  3393--3405.

\bibitem[BH]{Bridson-Haefliger}
M.~Bridson and A.~Haefliger, \emph{Metric spaces of non-positive curvature},
  Springer-Verlag, Berlin, 1999.

  \bibitem[Bri1]{Bridson:gwp}
  Martin~R. Bridson.
  \newblock {The geometry of the word problem}.
  \newblock In {\em Invitations to geometry and topology}, volume~7 of {\em Oxf.
    Grad. Texts Math.}, pages 29--91. Oxford Univ. Press, Oxford, 2002.

\bibitem[Bri2]{Bridson:PolynDehn}
M.~Bridson, \emph{Polynomial {D}ehn functions and the length of asynchronously
  automatic structures}, Proc. London Math. Soc. \textbf{85} (2002), no.~3,
  441--466.

  \bibitem[BT]{BurilloTaback:EquivalenceDehn}
  J.~Burillo and J.~Taback, \emph{Equivalence of geometric and
  combinatorial {D}ehn functions}, New York J. Math. \textbf{8} (2002),
  169--179.

\bibitem[DMS]{DrutuMozesSapir}
C.~Dru{\c{t}}u, Shahar Mozes, and Mark Sapir, \emph{Divergence in lattices in
  semisimple {L}ie groups and graphs of groups}, Trans. Amer. Math. Soc.
  \textbf{362} (2010), no.~5, 2451--2505.

\bibitem[ECH{\etalchar{+}}]{ECHLPT}
D.B.A. Epstein, J.~Cannon, D.F. Holt, S.~Levy, M.S. Paterson, and W.P.
  Thurston, \emph{{Word Processing and Group Theory}}, Jones and Bartlett,
  1992.

\bibitem[FF]{FedererFleming}
  H.~Federer and W.~Fleming, \emph{{Normal and integral currents}}, 
  Ann. Math.
    \textbf{72} (1960), 458--520.

\bibitem[Ger1]{Gersten:divergence3}
S.~Gersten, \emph{Divergence in 3-manifold groups}, Geom. Funct. Anal.
  \textbf{4} (1994), no.~6, 633--647.

\bibitem[Ger2]{Gersten:divergence}
\bysame, \emph{Quadratic divergence of geodesics in {CAT(0)}-spaces}, Geom.
  Funct. Anal. \textbf{4} (1994), no.~1, 37--51.

\bibitem[Gr1]{Gromov(1983)}
M.~Gromov, \emph{Filling {R}iemannian manifolds}, Journal of Differential
  Geometry \textbf{18} (1983), 1--147.

  \bibitem[Gr2]{Gromov:hyperbolic}
  M.~Gromov.
  \newblock {Hyperbolic groups}.
  \newblock In S.~Gersten, editor, {\em Essays in group theory}, volume~8 of {\em
    MSRI Publications}. Springer, 1987.

\bibitem[Gr3]{Gromov:Asymptotic}
\bysame, \emph{Asymptotic invariants of infinite groups}, Geometric Group
  Theory, Vol. 2 (Sussex, 1991) (G.~Niblo and M.~Roller, eds.), LMS Lecture
  Notes, vol. 182, Cambridge Univ. Press, 1993, pp.~1--295.

\bibitem[Grf1]{Groft:thesis}
Chad Groft, \emph{Isoperimetric functions on the universal covers of compact
  spaces}, ProQuest LLC, Ann Arbor, MI, 2007, Thesis (Ph.D.)--Stanford
  University.

\bibitem[Grf2]{Groft1}
\bysame, \emph{Generalized {D}ehn {F}unctions {I}}, preprint, \textsc{
  arXiv:0901.2303}.

\bibitem[Grf3]{Groft2}
\bysame, \emph{Generalized {D}ehn {F}unctions {II}}, preprint, \textsc{
  arXiv:0901.2317}.

\bibitem[Hat]{Hatcher:book}
A.~Hatcher, \emph{Algebraic {T}opology}, Cambridge University Press, 2002.


\bibitem[Hin]{Hindawi:div}
Mohamad~A. Hindawi, \emph{On the filling invariants at infinity of {H}adamard
  manifolds}, Geom. Dedicata \textbf{116} (2005), 67--85.

 \bibitem[Kle9]{Kleiner:lengthsp}
 Bruce Kleiner, \emph{The local structure of length spaces with curvature bounded above}, Math. Z. \textbf{231} (1999), no.~3, 409--456.

\bibitem[Lan]{Lang:linearHyp}
Urs Lang, \emph{Higher-dimensional linear isoperimetric inequalities in
  hyperbolic groups}, Internat. Math. Res. Notices (2000), no.~13, 709--717.

\bibitem[Leu1]{Leuzinger:divergence}
E.~Leuzinger, \emph{Corank and asymptotic filling-invariants for symmetric
  spaces}, Geom. Funct. Anal. \textbf{10} (2000), no.~4, 863--873.

\bibitem[Leu2]{Leuzinger:optimalDehn}
Enrico Leuzinger, \emph{Optimal higher-dimensional {D}ehn functions for some
  {CAT(0)} lattices}, preprint, \textsc{arXiv:1205.4923}, 2012.


\bibitem[Ols]{Olshanskii:hip}
 A.~Yu.~Olshanskii,\emph{ Hyperbolicity of groups with subquadratic isoperimetric inequalities}, Intl. J. Alg. Comp. \textbf{1} (1991), 282--290.


\bibitem[Pap95]{Papasoglu:hip}
P.~Papasoglu, \emph{On the subquadratic isoperimetric inequality}, Geometric group theory,
vol. \textbf{25}, de Gruyter, Berlin-New-York, 1995, R. Charney, M. Davis, M. Shapiro
(eds), pp. 193--200.

\bibitem[Pap96]{Papasoglu:quadratic}
P.~Papasoglu, \emph{On the asymptotic cone of groups satisfying a quadratic
  isoperimetric inequality}, Journal of Differential Geometry \textbf{44}
  (1996), 789--806.

\bibitem[Pap00]{Papasoglu:iso}
\bysame, \emph{Isodiametric and isoperimetric inequalities for complexes and
  groups}, J. London Math. Soc. (2) \textbf{62} (2000), no.~1, 97--106.

\bibitem[Ril]{TRiley:higherIso}
T.~R. Riley, \emph{Higher connectedness of asymptotic cones}, Topology
  \textbf{42} (2003), no.~6, 1289--1352.

\bibitem[Wen1]{Wenger:isoperim}
S.~Wenger, \emph{{Isoperimetric inequalities of Euclidean type in metric
  spaces}}, Geom. Funct. Anal. \textbf{15} (2005), 534--554.

\bibitem[Wen2]{Wenger:div}
\bysame, \emph{{Filling Invariants at Infinity and the Euclidean Rank of
  Hadamard Spaces}}, Int. Math. Res. Notices (2006), 1--33.

\bibitem[Wen3]{Wenger:asrk}
\bysame, \emph{The asymptotic rank of metric spaces}, Comment. Math. Helv.
  \textbf{86} (2011), no.~2, 247--275.

\bibitem[Wh]{White}
B.~White, \emph{Mappings that minimize area in their homotopy class}, J. Diff. Geom.
  \textbf{20} (1984), no.~2, 433--446.

\bibitem[You1]{young:nonrecursiveFV3}
R.~Young, \emph{Homological and homotopical higher order filling functions},
  preprint, \textsc{arXiv:0805.0584}.

  \bibitem[You2]{Young:slnz}
  R.~Young, \emph{The Dehn function of {SL(n;Z)}},
    Ann. of Math., to appear.

\end{thebibliography}
\end{document}